\documentclass{article}
\usepackage[utf8]{inputenc}
\usepackage{color,amsmath,bbm,mathtools,url}
\usepackage{xcolor}
\usepackage[normalem]{ulem}
\usepackage{bm}
\usepackage{layout}
\usepackage{subcaption}
\usepackage[capitalise,noabbrev]{cleveref}
\usepackage{cancel}
\crefformat{equation}{(#2#1#3)}
\numberwithin{equation}{section}
\numberwithin{figure}{section}
\usepackage{amsthm}
\usepackage{cancel}
\newcommand{\stkout}[1]{\ifmmode\text{\sout{\ensuremath{#1}}}\else\sout{#1}\fi}

\title{ A Time-Relaxation Reduced Order Model for the Turbulent Channel Flow}
\author{Ping-Hsuan Tsai, Paul Fischer, Traian Iliescu}
\date{\today}

\def\cN{{\mathcal{N}}}

\def\scriptO{{{\it O}\kern -.42em {\it `}\kern + .20em}}
\def\RR{{{\rm l}\kern - .15em {\rm R} }}
\def\PP{{{\rm l}\kern - .15em {\rm P} }}
\def\VV{{{\rm V}\kern - .69em {\rm V} }}
\def\L2{{{\sf L}^2}}
\def\H1{{{\sf H}^1}}
\def\PN2{{\PP_{N}-\PP_{N-2}}}
\def\QM2{{\mathcal{Q}_{M}-\mathcal{Q}_{M-2}}}

\def\complex{{{\rm C} \kern - .53em {\rm l} \kern + .38em}}

\def\a1{{ | \lambda_{\min} |}}

\def\l1{{   \lambda_{\min}  }}

\def\bff{{\bf f}}

\def\bu0{{\underline {\bf 0}}}

\def\bu{{\bm u}}

\def\bw{{\bf w}}
\def\bv{{\bf v}}
\def\bx{{\bf x}}

\def\bX{{\bf X}}

\def\bn{{\bf n}}

\def\uf{{\underline f}}

\def\uu{{\underline u}}

\def\u0{{\underline 0}}

\newcommand{\pp}[2]{\frac{\partial #1}{\partial #2} }

\def\bX{{\bf X}}

\def\ubu{{\overline{\underline u}}}

\def\ubu{{\underline{\overline{u}}}}

\def\tr{\text{r}}

\def\tProj{\text{Proj}}
\def\trecon{\text{recon}}
\def\tpred{\text{pred}}

\def\bracu{\langle u \rangle}
\def\bracur{\langle u_\tr \rangle}
\def\bracvr{\langle v_\tr \rangle}

\def\bracuproj{\langle u_\tProj \rangle}
\def\bracvproj{\langle v_\tProj \rangle}
\def\bracv{\langle v \rangle}

\def\ou{\overline{u}}
\def\obu{\overline{\bu}}

 \def\euu{\varepsilon_{u'u'}} 
 \def\euv{\varepsilon_{u'v'}}

\newtheorem{remark}{Remark}[section]

\def\PP{{{\rm l}\kern - .15em {\rm P} }}
\def\PN2{{\PP_{N}-\PP_{N-2}}}

\newcommand{\cO}{\mathcal{O}}

\newcommand{\bphi}{\boldsymbol{\varphi}}

\definecolor{vargreen}{rgb}{0.0, 0.5, 0.0}

\newcommand{\deleted}[1]{{}}

\begin{document}

\maketitle
\begin{abstract}
Regularized reduced order models (Reg-ROMs) are stabilization strategies that
leverage spatial filtering to alleviate the spurious numerical oscillations
generally displayed by the classical Galerkin ROM (G-ROM) in under-resolved
numerical simulations of turbulent flows.  
In this paper, we propose a new Reg-ROM, the time-relaxation ROM (TR-ROM), which
filters the marginally resolved scales.  
We compare the new TR-ROM with the two other Reg-ROMs in current use, i.e., the
Leray ROM (L-ROM) and the evolve-filter-relax ROM (EFR-ROM), in the numerical
simulation of the turbulent channel flow at $Re_{\tau} = 180$ and $Re_{\tau} =
395$ in both the reproduction and the predictive regimes.
For each Reg-ROM, we investigate two different filters:
    (i) the differential filter (DF), and
    (ii) a new higher-order algebraic filter (HOAF). 
In our numerical investigation, we monitor the Reg-ROM performance with respect
to the ROM dimension, $N$, and the filter order.  We also perform sensitivity
studies of the three Reg-ROMs with respect to the time interval, relaxation
parameter, and filter radius.  The numerical results yield the following
conclusions:
    (i) All three Reg-ROMs are significantly more accurate than the G-ROM. 
    (ii) All three Reg-ROMs are more accurate than the ROM projection, which
    represents the best theoretical approximation of the training data in the
    given ROM space. 
    (iii) With the optimal parameter 
    values, the new
    TR-ROM yields more accurate results than the L-ROM and the EFR-ROM in all tests.
    (iv) For most $N$ values, DF yields the most accurate results for all three
    Reg-ROMs.  
    (v) The optimal parameters trained in the reproduction regime are also
    optimal for the predictive regime for most $N$ values, demonstrating the
    Reg-ROM predictive capabilities.  
    (vi) All three Reg-ROMs are sensitive to the filter radius and the filter
    order, and the EFR-ROM and the TR-ROM are sensitive to the relaxation
    parameter.  (vii) The optimal range for the filter radius and the effect of
    relaxation parameter are similar for the two $\rm Re_\tau$ values.
\end{abstract}


\section{Introduction}

Reduced order models (ROMs) are computational models that leverage data to
approximate the dynamics in a space whose dimension is orders of magnitude lower
than the dimension of full order models (FOMs), i.e., models obtained from
classical numerical discretizations (e.g., the finite element or spectral
element methods).  In the numerical simulation of fluid flows, Galerkin ROMs
(G-ROMs), which use a data-driven basis in a Galerkin framework, have provided
efficient and accurate approximations of laminar flows, such as the
two-dimensional flow past a circular cylinder at low Reynolds
numbers~\cite{brunton2019data,hesthaven2015certified,noack2011reduced,quarteroni2015reduced}.

However, turbulent flows (e.g., the turbulent channel
flow~\cite{kim1987turbulence,moser1999direct}) are notoriously hard for the
standard G-ROM: To capture the complex dynamics of the turbulent flow, a
relatively large number (on the order of hundreds and even thousands~\cite[Table
II]{ahmed2021closures}\cite{tsai2023accelerating}) of ROM basis functions are needed.  Thus, the resulting
G-ROM is relatively high-dimensional and its computational cost is too large to
be used in realistic applications, such as, control of turbulent flows.  To
mitigate this issue, the standard G-ROM is often constructed with relatively few
basis functions.  The resulting G-ROM is appealing because it is low-dimensional
and computationally efficient.  It is, however, inaccurate.  The reason is that
the ROM basis functions that were not used to build the low-dimensional G-ROM
have an important role in the G-ROM dynamics.  Indeed, as shown in~\cite{CSB03},
the role of the discarded ROM modes is to extract energy from the system.
Without including this dissipation mechanism, the under-resolved G-ROM (i.e.,
the G-ROM that does not include enough basis functions to capture the underlying
complex dynamics) generally yields spurious numerical oscillations.  Thus, in
the numerical simulation of turbulent flows, these efficient, low-dimensional
ROMs are generally equipped with ROM closures (see the review
in~\cite{ahmed2021closures}) and stabilizations (see,
e.g.,~\cite{carlberg2011efficient,parish2020adjoint,mou2023energy,kaneko2022augmented,fick2018stabilized}).
For example, in the reduced order modeling of the turbulent channel flow, ROM
closures (e.g., the eddy viscosity ROM~\cite{johansson2006reduced} and the
mixing-length ROM~\cite{mou2023energy}) or ROM stabilizations (e.g., the
evolve-filter-relax ROM~\cite{mou2023energy}) have been successfully employed.
Regularized ROMs (Reg-ROMs) are stabilizations that leverage ROM spatial
filtering of terms of the Navier-Stokes equations to decrease the size of the
spurious numerical oscillations and increase the ROM accuracy.  The two Reg-ROMs
in current use are (i) the Leray ROM (L-ROM)
\cite{wells2017evolve,sabetghadam2012alpha,strazzullo2022consistency,kaneko2020towards},
in which the velocity component of the convective term of the Navier-Stokes
equations is filtered, and (ii) the evolve-filter-relax ROM (EFR-ROM)
\cite{wells2017evolve,strazzullo2022consistency}, which filters an intermediate
approximation of the Navier-Stokes equations, and then relaxes it.

In this paper, we propose a new Reg-ROM, the time-relaxation ROM (TR-ROM), which
filters the marginally resolved scales.  The main goal of this paper is to
investigate the new TR-ROM in the numerical simulation of the turbulent channel
flow at friction Reynolds numbers $\rm Re_{\tau} = 180$ and $395$, in both the
reproduction and the predictive regimes.  To ensure a thorough assessment of the
new TR-ROM, we compare it with the two other Reg-ROMs in current use, i.e.,
L-ROM and EFR-ROM.  For each Reg-ROM, we investigate two different filters: (i)
the differential filter (DF); and (ii) a new higher-order algebraic filter.  
We also compare the three Reg-ROMs with the classical G-ROM
(i.e., the ROM that does not use any stabilization or closure).  As a benchmark for
the three Reg-ROMs and the G-ROM, we use the FOM, which is a direct numerical
simulation (DNS) of the turbulent channel flow.  Furthermore, in our numerical
comparison, we use the projection error, which is the error with respect to the
projection of the FOM solution onto the space spanned by the ROM basis
functions.  
We expect the three Reg-ROM to alleviate the G-ROM spurious numerical
oscillations and yield more accurate results with a negligible computational
overhead. 

It is well known that, in numerical simulation of turbulent flows, the
parameters used in the computational setting can have a significant effect on
the ROM results.  Thus, to ensure a fair assessment of the Reg-ROMs, we perform
a sensitivity study of the numerical results with respect to the following
parameters:
(i) the number of ROM basis functions, $N$, used to construct the ROM; 
(ii) the higher-order algebraic filter order, $m$;
(iii) the filter radius, $\delta$;
(iv) the relaxation parameter, $\chi$; 
and (v) the time interval.

The rest of the paper is organized as follows: 
In Section~\ref{section:numerical-modeling}, we present the FOM and the G-ROM.
In Section~\ref{section:rom-filters}, we describe the two ROM spatial filters
that we use to construct the Reg-ROMs: the classical differential filter and the
higher-order algebraic filter.  
In Section~\ref{section:reg-roms}, we leverage these two ROM filters to
construct the L-ROM, EFR-ROM, and the new TR-ROM.  
In Section~\ref{section:numerical-results}, we perform a numerical investigation
of the three Reg-ROMs (i.e., L-ROM, EFR-ROM, and TR-ROM) in the numerical
simulation of the turbulent channel flow at friction Reynolds numbers $\rm
Re_{\tau}=180$ and $395$. 
In addition, we assess the accuracy of the three Reg-ROMs and their
sensitivity with respect to parameters in both the reproduction and the
predictive regimes.  
In Section~\ref{section:conclusions}, we present the
conclusions of our numerical investigation and we outline directions of future
research.  
Finally, in Appendix~\ref{section:appendix}, we present a theoretical
and numerical investigation of the higher-order algebraic filter.


\section{Numerical Modeling}
    \label{section:numerical-modeling}

In this section, we briefly outline the FOM (Section~\ref{section:fom}) and the
G-ROM (Section~\ref{section:g-rom}) used in our numerical investigation.

\subsection{Full Order Model (FOM)}
    \label{section:fom}
    
The governing equations are the incompressible Navier-Stokes equations (NSE) with forcing: 
\begin{align} \label{eq:PDE}
  \frac{\partial \bu}{\partial t} + (\bu \cdot \nabla) \bu &=  - \nabla p + \frac{1}{\rm Re}
   \nabla^2 \bu + \bff(\bu), \qquad \nabla \cdot \bu = 0, 
\end{align} 
where $\bu$ is the velocity and $p$ the pressure.  Here, $\bff(\bu)$ is a
uniform forcing vector field function in the streamwise direction, $x$, that
enforces a time-constant flow-rate on the solution. The initial condition
consists of random noise, which eventually triggers transition to turbulence,
and the boundary conditions are periodic in the streamwise and spanwise
directions of the channel, and homogeneous Dirichlet in the wall-normal
direction.

The FOM is constructed using the Galerkin projection of \cref{eq:PDE} onto
the spectral element space with the $\mathcal{P}_q-\mathcal{P}_q$
velocity-pressure coupling. Following \cite{fischer2017recent}, a semi-implicit
scheme BDF$k$/EXT$k$ is used for time discretization.  Specifically, the
$k$th-order backward differencing (BDF$k$) is used for the time-derivative term,
$k$th-order extrapolation (EXT$k$) for the advection and forcing terms, and
implicit treatment on the dissipation terms. As discussed in
\cite{fischer2017recent}, $k=3$ is used to ensure the imaginary eigenvalues
associated with skew-symmetric advection operator are  within the stability
region of the BDF$k$/EXT$k$ time-stepper. 

The full discretization leads to solving a linear unsteady Stokes system at each
time step. The forcing term effectively adds an impulse-response streamwise
velocity field.  This impulse response is scaled appropriately at each time step
to ensure that the mean velocity at each timestep yields the prescribed
flow-rate \cite{ghaddar1986conservative,patankar1977fully}.  The detailed
derivation of the FOM and the treatment of the constant-flow rate can be found
in \cite{tsai2022parametric} and \cite{kaneko2022augmented}, respectively. 

\subsection{Galerkin Reduced Order Model (G-ROM)}
    \label{section:g-rom}
    
In this section, we introduce the G-ROM.  We follow the standard proper
orthogonal decomposition (POD) procedure
\cite{berkooz1993proper,volkwein2013proper} to construct the reduced basis function.  To
this end, we collect a set of DNS solutions lifted by the zeroth mode $\bphi_0$,
and form its corresponding Gramian matrix using the $L^2$ inner product (see,
e.g.,~\cite{fick2018stabilized,kaneko2020towards} for alternative strategies).
The first $N$ POD basis functions $\{\bphi_i\}^N_{i=1}$ are constructed from the first $N$
eigenmodes of the Gramian.  Setting the zeroth mode, $\bphi_0$, to the
time-averaged velocity field in the time interval in which the snapshots were
collected, the G-ROM is constructed by inserting the ROM basis expansion
\begin{equation} \label{eq:romu}
   \bu_\tr(\bx) = \bphi_0(\bx) + \sum_{j=1}^N u_{\tr,j} \bphi_j(\bx)
\end{equation}
into the weak form of \cref{eq:PDE}:
{\em Find $\bu_\tr$ 
such that, for all $\bv \in \bX^N_0$,}
\begin{eqnarray}
    && 
    \left(
        \frac{\partial \bu_\tr}{\partial t} , \bv_{i} 
    \right)
    + Re^{-1} \, 
    \left( 
        \nabla \bu_\tr , 
        \nabla \bv_{i} 
    \right)
    + \biggl( 
        (\bu_\tr \cdot \nabla) \bu_\tr ,
        \bv_{i} 
    \biggr)
    = 0,  
    \label{eq:gromu}
\end{eqnarray}
where $(\cdot,\cdot)$ denotes the $L^2$ inner product and 
$\bX^N_0 := \text{span} \{\bphi_i\}^N_{i=1}$ is the ROM space. 

\begin{remark}
    We note that, in the case of fixed geometries, the divergence and pressure
    terms drop out of \cref{eq:gromu} because the ROM basis function is weakly
    divergence-free. For ROMs that include the pressure approximation, see,
    e.g.,~\cite{hesthaven2015certified,quarteroni2015reduced,ballarin2015supremizer,decaria2020artificial,noack2005need}.
\end{remark}

\begin{remark}
    Because the zeroth mode has the prescribed flow-rate, the remaining POD
    basis functions have zero flow-rate, meaning that the test space $\bX^N_0 =
    \text{span}\{\bphi\}^N_{i=1}$ contains only members with zero flow-rate.
    Hence no additional forcing term is required for the ROM formulation because
    the forcing term drops out of (\ref{eq:gromu}):
    \begin{equation}
        \int_{\Omega} \bv \cdot \bff\, dV = f_x \left( \int_{\Omega} v_x dV\right) = 0,\quad \forall~ \bv \in \bX^N_0, 
    \end{equation}
    where $f_x$ and $v_x$ are the streamwise components of $\bff$ and $\bv$, respectively.
\end{remark}

With (\ref{eq:gromu}), the following evolution equations are
derived for the ROM basis coefficients $u_{\tr,j}$: For each $i=1,\ldots,N$,
\begin{align}
	\sum^N_{j=1}B_{ij}\frac{d u_{\tr,j}}{dt} & =
	-\sum^N_{k=0}\sum^N_{j=0} C_{ikj} u_{\tr,k}(t) u_{\tr,j}(t) - \frac{1}{\rm Re}
	\sum^N_{j=0} A_{ij} u_{\tr,j}(t), \label{eq:nse_ode}
\end{align}
where $A$, $B$, and $C$ represent the stiffness, mass, and advection operators,
respectively, with entries
\begin{align} 
    A_{ij} &=\int_{\Omega}\nabla \bphi_i :\nabla \bphi_j\, dV, \quad  B_{ij} =\int_{\Omega}\bphi_i\cdot\bphi_j\, dV, \label{eq:Bu} \\ 
    C_{ikj}& =\int_{\Omega}\bphi_i\cdot(\bphi_k\cdot\nabla)\bphi_j\,dV. \label{eq:Cu} 
\end{align}

We note that, because the POD basis functions are orthonormal in the
$L^2$ inner product, the ROM mass matrix $B$ is an identity matrix.

Applying the BDF$k$/EXT$k$ scheme to the reduced system (\ref{eq:nse_ode}), the
fully discrete $N \times N$ reduced system can be written as follows:
\begin{align}
	H_{\tr,\rm Re} \, \uu_{\tr}^{l+1} & =
	{\uf}(\uu_{\tr}^l;\rm Re)\label{eq:nse_d1},
\end{align}
for each timestep $t^l$, where $H_{\tr, \rm Re}$ and $\uf(\uu_{\tr}^l;\rm Re)$
are the resulting Helmholtz matrix and right-hand side vector, and $\uu_{\tr}^l$
is the vector of ROM coefficients of $\bu_{\tr}^l$.

\section{ROM Filters}
    \label{section:rom-filters}
In this section, we present the two ROM spatial filters that we use to construct
the Reg-ROMs in Section~\ref{section:reg-roms}: the classical differential
filter (Section~\ref{section:df}) and the higher-order algebraic filter
(Section~\ref{section:hodf}).

\subsection{ROM Differential Filter (DF)}
    \label{section:df}
The first ROM spatial filter we investigate is the {\it ROM differential filter (DF)}:
{\em Given $\bu_\tr = \sum^N_{j=1} u_{\tr,j} \bphi_j(\bx)$, find 
$\obu_\tr(\bx) = \sum_{j=1}^N \ou_{\tr,j} \bphi_j(\bx)$
such that}
\begin{eqnarray} 
    \biggl( 
        \obu_\tr - \delta^2 \Delta 
        \obu_\tr , \bphi_i 
        \biggr) 
        = \biggl(\bu_\tr, \bphi_i \biggr)
    \qquad \forall \, i=1, \ldots N,
    \label{equation:df-weak}
\end{eqnarray}
where $\delta$ is the filter radius.  The DF weak form~\eqref{equation:df-weak}
yields the following linear system:
\begin{eqnarray} 
    \left( \mathbbm{I} + \delta^2 
    A \right) \underline{\ou}_{\tr}
    = \underline{u}_{\tr},
    \label{equation:df-linear-system}
\end{eqnarray}
where $\underline{\ou}_{\tr}$ is the vector of ROM coefficients of $\obu_\tr$,
and $\mathbbm{I}$ and $A$ are the identity and ROM stiffness matrices,
respectively.  \footnote{In general, the inverse of the ROM mass matrix,
$B^{-1}$, shows up in the DF definition.  However, because the POD basis function is
orthonormal in the $L^2$ norm, $B$ is the identity matrix and ignored here to
avoid confusion.} 
We note that the expansions for $\bu_\tr$ and
$\obu_\tr$ do not include the zeroth mode, $\bphi_{0}$.  This is in contrast
with the expansion~\eqref{eq:romu}, which does include $\bphi_{0}$.  The reason
for not including $\bphi_{0}$ in our expansions is that 
this strategy
was shown in~\cite{wells2017evolve} to yield more accurate results. 

We emphasize that~\eqref{equation:df-linear-system} is a low-dimensional, $N
\times N$ linear system, whose computational overhead is negligible.  Thus, DF
will be used in Section~\ref{section:reg-roms} to construct regularized ROMs
that increase the ROM accuracy without significantly increasing the
computational cost.  

DF was used in large eddy simulation of turbulent flows with classical numerical
discretizations~\cite{germano1986differential,layton2012approximate}.  In
reduced order modeling, DF was used to develop Reg-ROMs for the
Kuramoto-Sivashinsky equation~\cite{sabetghadam2012alpha}, the
NSE~\cite{wells2017evolve,strazzullo2022consistency}, and the quasi-geostrophich
equations~\cite{girfoglio2023linear}.

\subsection{ROM Higher-Order Algebraic Filter (HOAF)} 
    \label{section:hodf}

The second filter we investigate is the {\it higher-order algebraic filter (HOAF)}: 
{\em Given $\bu_\tr = \sum^N_{j=1} u_{\tr,j} \bphi_j(\bx)$, find $\obu_\tr(\bx)
= \sum_{j=1}^N \ou_{\tr,j} \bphi_j(\bx)$
such that}
\begin{eqnarray} 
    H_{m} \underline{\ou}_{\tr}
    \coloneqq \left( \mathbbm{I} + \delta^{2m} A^{m} \right) \underline{\ou}_{\tr}
    = \uu_{\tr}, 	
    \label{eqn:hodf}
\end{eqnarray}
where $\underline{\ou}_{\tr}$ is the vector of ROM basis coefficients of
$\obu_\tr$, and $m$ is a positive integer.  
As motivated in Section~\ref{section:df}, the expansions for $\bu_\tr$ and
$\obu_\tr$ do not include the zeroth mode, $\bphi_{0}$.  As explained
in~\cite{mullen1999filtering} in the Fourier setting, the role of the exponent
$m$ is to control the percentage of filtering at different wavenumbers: As $m$
increases, the amount of filtering increases for the high wavenumber components
and decreases for the low wavenumber components.  

Just as the DF~\eqref{equation:df-linear-system}, the HOAF~\eqref{eqn:hodf} is
also a low-dimensional, $N \times N$ linear system.  Thus, HOAF will also be
used in Section~\ref{section:reg-roms} to develop accurate and efficient
Reg-ROMs.

\begin{remark}[Notation Convention]
    We also note that, for $m=1$, the linear systems~\eqref{eqn:hodf} and
    \eqref{equation:df-linear-system} are identical.  Thus, DF can be considered
    a particular case of HOAF with $m=1$.  In what follows, for notation
    convenience, we will use the linear system~\eqref{eqn:hodf} for both HOAF
    and DF, and we will only specify the $m$ value to differentiate between the
    two: $m=1$ for DF and $m\geq2$ for HOAF.  \label{remark:hodf-m1} 
\end{remark}

\begin{remark}
    The HOAF~\eqref{eqn:hodf} was proposed in~\cite{gunzburger2019evolve} and
    was based on the HOAF introduced by Fischer and
    Mullen~\cite{mullen1999filtering} in a spectral element method (SEM)
    setting.  We also note that the HOAF used in~\cite{gunzburger2019evolve} has
    a $\delta$ scaling that is dimensionally inconsistent.  This is rectified
    in~\eqref{eqn:hodf}. 
\end{remark}

\begin{remark}[Nomenclature]
    In~\cite{gunzburger2019evolve}, the HOAF~\eqref{eqn:hodf} was called the
    higher-order differential filter, to be consistent with the SEM
    nomenclature.  In Section~\ref{section:theoretical-investigation}, we show
    that the HOAF is related to, but slightly different from, the spatial
    discretization of a higher-order differential operator.  
    (They are the same
    in the periodic case.) Thus, for clarity, 
    in this paper we call the operator in~\eqref{eqn:hodf} the high-order
    algebraic filter: The term $A^{m}$ yields the high-order algebraic character
    of the operator~\eqref{eqn:hodf}, and the numerical investigation in
    Section~\ref{section:numerical-investigation} shows that the
    operator~\eqref{eqn:hodf} acts like a spatial filter. 
\end{remark}

A theoretical and numerical investigation of HOAF is performed in Appendix~\ref{section:appendix}.

\section{Regularized Reduced Order Models (Reg-ROMs)}
    \label{section:reg-roms}
    
In this section, we outline the three Reg-ROMs that we compare in our numerical
investigation in Section~\ref{section:numerical-results}: the Leray ROM
(Section~\ref{section:l-rom}), the evolve-filter-relax ROM
(Section~\ref{section:efr-rom}), and the novel time relaxation ROM
(Section~\ref{section:tr-rom}).  All three Reg-ROMs are developed based on the
same principle: Use the ROM spatial filters presented in
Section~\ref{section:rom-filters} to smooth (filter) terms in the standard
G-ROM~\eqref{eq:gromu} and eliminate/alleviate the G-ROM's spurious numerical
oscillations in under-resolved turbulent flow simulations.  Furthermore, because
the computational cost of the ROM spatial filters is low, the Reg-ROM
computational overhead with respect to the G-ROM is negligible.  Thus, Reg-ROMs
are expected to yield more accurate results than the standard G-ROM without
a significantly increase in the computational cost.

\subsection{Leray ROM (L-ROM)}
    \label{section:l-rom}
    
The {\it Leray ROM (L-ROM)}~\cite{kaneko2020towards,wells2017evolve} modifies
the standard G-ROM weak formulation~\eqref{eq:gromu} as follows: {\em Find
$\bu_{\tr}$ of the form~\eqref{eq:romu} such that, $\forall \, i=1, \ldots N,$}
\begin{eqnarray}
    && 
    \left(
        \frac{\partial \bu_\tr}{\partial t} , \bphi_{i} 
    \right)
    + Re^{-1} \, 
    \left( 
        \nabla \bu_\tr , 
        \nabla \bphi_{i} 
    \right)
    + \biggl( 
        (\obu_\tr \cdot \nabla) \bu_\tr ,
        \bphi_{i} 
    \biggr)
    = 0  
    \label{eqn:l-rom}
\end{eqnarray}
where $\obu_\tr$ is the ROM velocity filtered with one of the ROM spatial
filters introduced in Section~\ref{section:rom-filters}, that is,
DF~\eqref{equation:df-linear-system} or HOAF~\eqref{eqn:hodf}.

L-ROM~\eqref{eqn:l-rom} is a Reg-ROM, because it leverages spatial filtering of
the convective term of the G-ROM (i.e., it replaces $(\bu_\tr \cdot \nabla)
\bu_\tr$ with $(\obu_\tr \cdot \nabla) \bu_\tr$) in order to smooth out the
G-ROM's spurious numerical oscillations in the convection-dominated,
under-resolved regime. 

In a more general setting, the Leray model was first
introduced by Jean Leray in 1934 as a theoretical tool 
to prove existence of weak solutions of the NSE~\cite{leray1934sur}.
As a computational tool, Leray regularization was first used
in~\cite{geurts2003regularization} as a stabilization strategy for
under-resolved simulations of turbulent flows with classical numerical
discretizations~\cite{layton2012approximate}.  
As noted by Guermond and co-authors 
\cite{guermond2004mathematical,guermond2011entropy},
when a differential filter is used, the Leray model is similar to the
NS-$\alpha$ model of Foias, Holm, and Titi~\cite{FHT2001}.
Leray regularization was first used in the context of
reduced order models in~\cite{sabetghadam2012alpha} 
for the Kuramoto-Sivashinsky equations.  
For fluid flows, L-ROM was first used in~\cite{wells2017evolve} for the 3D flow
past a circular cylinder at $Re=1000$.  Since then, L-ROM has been
successfully used as a stabilization technique for various under-resolved
flows: the NSE~\cite{girfoglio2021pod,girfoglio2023hybrid}, the stochastic
NSE~\cite{gunzburger2019evolve,gunzburger2020leray}, and the quasigeostrophic
equations~\cite{girfoglio2023linear,girfoglio2023novel}.  
To our knowledge, L-ROM has never been used for 
the turbulent channel flow.

\subsection{Evolve-Filter-Relax ROM (EFR-ROM)}
    \label{section:efr-rom}

The {\it evolve-filter-relax ROM (EFR-ROM)}, introduced in
~\cite{wells2017evolve,gunzburger2019evolve}, consists of three steps.
{\em Given the EFR-ROM
approximation at the current time step, $\bu_{\tr}^{n}$, find the EFR-ROM
approximation at the next time step, $\bu_{\tr}^{n+1}$, as follows:}
\begin{eqnarray}
   && \hspace*{-1.2cm} 
   \text{\bf (I)} \text{\emph{ Evolve}:} \quad \text{{\em Find $\bw_{\tr}$ of the form~\eqref{eq:romu} such that, $\forall \, i=1, \ldots N,$}} \nonumber \\[0.4cm]	
        && \hspace*{-0.5cm} \left(
            \frac{\bw_{\tr}^{n+1} - \bu_{\tr}^{n}}{\Delta t} , \bphi_{i} 
        \right)
        + Re^{-1} \, 
        \left( 
            \nabla \bu_{\tr}^{n} , 
            \nabla \bphi_{i} 
        \right)
        + \biggl( 
            (\bu_{\tr}^{n} \cdot \nabla) \bu_{\tr}^{n} , \bphi_{i} 
        \biggr)
        = 0  \\[0.4cm]	
 && \hspace*{-1.2cm} \text{\bf (II)} \text{\emph{ Filter}:} \quad 
    \bw_{\tr}^{n+1} 
    \longmapsto \overline{\bw_{\tr}^{n+1}}  \\[0.6cm]
     && \hspace*{-1.2cm} 
    \text{\bf (III)} \text{\emph{ Relax}:} \quad 
    \bu_{\tr}^{n+1} 
    = (1 - \chi) \, \bw_{\tr}^{n+1}
            + \chi \, \overline{\bw}_{\tr}^{n+1}. \label{eqn:efr-rom}
\end{eqnarray}

    In Step (I) of the EFR-ROM, called the evolve step, one step of the standard G-ROM
    time discretization is used to advance the current EFR-ROM approximation,
    $\bu_{\tr}^{n}$, to an intermediate EFR-ROM approximation, $\bw_{\tr}^{n+1}$. 
    In Step (II), called the filter step, one of the two ROM spatial
    filters presented in Section~\ref{section:rom-filters} is used to filter the
    intermediate EFR-ROM approximation obtained in Step (I) and obtain a smoother
    approximation, without spurious numerical oscillations.  Finally, in Step (III),
    called the relax step, the EFR-ROM approximation at the next time
    step is defined as a convex combination of the unfiltered intermediate
    EFR-ROM approximation obtained in Step (I), $\bw_{\tr}^{n+1}$, and its filtered
    counterpart, $\overline{\bw}_{\tr}^{n+1}$.  The goal of the relax step is to adjust
    the amount of dissipation introduced in the filter step by using a relaxation
    parameter, $0 \leq \chi \leq 1$. By varying $\chi$, one can produce a
    range of filter strengths, from no filtering at all ($\chi=0$) to 
    maximum filtering ($\chi=1$).  
    We note that the numerical investigation in~\cite{strazzullo2022consistency}
    has shown that EFR-ROM is sensitive with respect to $\chi$.

    The EFR strategy is well developed for classical numerical discretizations, 
    for example, in the context of
    finite element method~\cite{layton2012approximate}, and the spectral and
    spectral element methods~\cite{mullen1999filtering,DFM02}.  In reduced order
    modeling, the evolve-filter ROM was introduced in~\cite{wells2017evolve} and
    EFR-ROM in~\cite{gunzburger2019evolve}.  Since then, EFR-ROM has been developed
    in several directions, for example, the FOM-ROM
    consistency~\cite{strazzullo2022consistency} and feedback
    control~\cite{strazzullo2023new}.

\subsection{Time Relaxation ROM (TR-ROM)}
    \label{section:tr-rom}

In this paper, we propose a new type of Reg-ROM: the {\it time-relaxation ROM
(TR-ROM)}: {\em Find $\bu_{\tr}$ of the form~\eqref{eq:romu} such that, $\forall
\, i=1, \ldots N,$}
\begin{align}
    \left( \pp{\bu_\tr}{t} , \bphi_{i} \right)
    + \frac{1}{\rm Re} \left( \nabla \bu_\tr , \nabla \bphi_{i} \right)
    + \biggl( (\bu_\tr \cdot \nabla) \bu_\tr , \bphi_{i}  \biggr)
    + \biggl( \chi (\bu_\tr - 
    \obu_\tr) , \bphi_{i}  \biggr)
    = 0 ,
    \label{eqn:tr-rom}
\end{align}
where $\chi$ is the time-relaxation parameter, and $\obu_\tr$ is the ROM
velocity filtered with one of the ROM spatial filters introduced in
Section~\ref{section:rom-filters}, that is, DF~\eqref{equation:df-linear-system}
or HOAF~\eqref{eqn:hodf}.

Time-relaxation has been used as a regularization/stabilization strategy at a
FOM level~\cite[Chapter 5]{layton2012approximate} (see
also~\cite{SAK01a,SAK01b,AS02,layton2007truncation,schlatter2005evaluation}).
To our knowledge, this is the first use of time-relaxation stabilization in the
ROM context.  

To understand the role of the time relaxation term, we consider $\bu_\tr$ as
test function in TR-ROM~\eqref{eqn:tr-rom}, which is the usual approach in
energy-stability analysis. With this
choice, the last term on the right-hand side of~\eqref{eqn:tr-rom} can be
written as follows:
\begin{eqnarray}
    ( \chi (\bu_\tr - \obu_\tr) , \bu_\tr )
    = \chi ( \bu'_\tr , \bu_\tr ),
    \label{eqn:tr-rom-1}
\end{eqnarray}
where $\bu'_\tr$ represents the fluctuations of $\bu_\tr$ around $\obu_\tr$:
\begin{eqnarray}
    \bu'_\tr
    := \bu_\tr - \obu_\tr.
    \label{eqn:tr-rom-2}
\end{eqnarray}
Using the decomposition~\eqref{eqn:tr-rom-2}, the inner product in the time relaxation term~\eqref{eqn:tr-rom-1} becomes:
\begin{eqnarray}
    ( \bu'_\tr , \bu_\tr )
    = ( \bu'_\tr , \bu'_\tr )
    + ( \bu'_\tr , \obu_\tr ).
    \label{eqn:tr-rom-3}
\end{eqnarray}
The first term on the right-hand side of~\eqref{eqn:tr-rom-3} has a clear physical interpretation:
It is a dissipative term acting only on the fluctuations.  To understand the
role played by the second term on the right-hand side of~\eqref{eqn:tr-rom-3},
we distinguish two cases:

{\em Case 1: DF} \ 
When DF is used to construct TR-ROM, we can use~\eqref{equation:df-weak} and \eqref{eqn:tr-rom-2} to formally write the following: 
\begin{eqnarray}
    \bu'_\tr
    = \bu_\tr - \obu_\tr
    = - \delta^2 \Delta \obu_\tr.
    \label{eqn:tr-rom-4}
\end{eqnarray}
Using~\eqref{eqn:tr-rom-4} in the last term in~\eqref{eqn:tr-rom-3}, we obtain the following:
\begin{eqnarray}
    ( \bu'_\tr , \bu_\tr )
    &=& ( \bu'_\tr , \bu'_\tr )
    + ( \bu'_\tr , \obu_\tr )
    \nonumber \\
    &=&  ( \bu'_\tr , \bu'_\tr )
    - \delta^2 ( \Delta \obu_\tr , \obu_\tr )
    \nonumber \\
    &\stackrel{\obu_\tr = 0 \text{ on } \partial \Omega}{=}&  ( \bu'_\tr , \bu'_\tr )
    + \delta^2 ( \nabla \obu_\tr , \nabla \obu_\tr )
    \nonumber \\
    &\stackrel{\eqref{eqn:tr-rom-4}}{=}&  \delta^4 ( \Delta \obu_\tr , \Delta \obu_\tr )
    + \delta^2 ( \nabla \obu_\tr , \nabla \obu_\tr ).
    \label{eqn:tr-rom-5}
\end{eqnarray}
Equality~\eqref{eqn:tr-rom-5} shows that the TR-ROM term in~\eqref{eqn:tr-rom-1}
is a dissipative term that has two components: The first component can be
interpreted either as a dissipation term acting on the fluctuations or as a
hyperviscosity term acting on the averages.  The second component is a diffusion
term acting on the averages. 

\bigskip

{\em Case 2: HOAF} \ 
The HOAF is a high-order algebraic filter (\ref{eqn:hodf}) that does not have a
direct interpretation as a differential operator (primarily because of ambiguity
related to boundary conditions; see discussion in Appendix~\ref{section:appendix}).
When HOAF (\ref{eqn:hodf}) is used to construct TR-ROM, we cannot use the above
approach to interpret the TR-ROM term because the HOAF cannot be easily written in
terms of the spatial derivative operators. 
However, we can (optimistically) expect behavior of the following form
for a $2m$th-order filter,
\begin{eqnarray}
    ( \bu'_\tr , \bu_\tr )
    &=& ( \bu'_\tr , \bu'_\tr )
    + \delta^{2m} ( \nabla^{m} \obu_\tr , \nabla^{m} \obu_\tr )
    \nonumber \\
    &=& \delta^{4m} ( \Delta^{m} \obu_\tr , \Delta^{m} \obu_\tr )
    + \delta^{2m} ( \nabla^{m} \obu_\tr , \nabla^{m} \obu_\tr ),
    \label{eqn:tr-rom-6}
\end{eqnarray}
which has the same physical interpretation as that in the DF case 
and which is in fact exact if the domain is periodic in all directions.

\begin{remark}
    \label{remark:ad-vs-hodf}
We note that, at the FOM level, various levels of spatial filtering in the time
relaxation term has been obtained by using approximate deconvolution strategies
of different orders~\cite[Chapter 5]{layton2012approximate}.  In the new TR-ROM,
we employ a different strategy and adjust the HOAF order in order to control the
amount of filtering in TR-ROM.  Although the zeroth order deconvolution and the
DF yield identical models~\cite{layton2012approximate}, the higher-order
approximate deconvolution methods and the HOAF yield different time relaxation
models.

A numerical comparison of DF and HOAF, and an investigation of their effect on
the TR-ROM results is performed in Section~\ref{section:numerical-results}.
\end{remark}

\section{Numerical Results}
    \label{section:numerical-results}
    
In this section, we present numerical results for the three Reg-ROMs outlined in
Section~\ref{section:reg-roms} in the simulation of the turbulent channel flow.
Specifically, we compare the L-ROM~\eqref{eqn:l-rom}, the
EFR-ROM~\eqref{eqn:efr-rom}, and the novel TR-ROM~\eqref{eqn:tr-rom}.
For comparison purposes, we also investigate the standard G-ROM~\eqref{eq:gromu}.
As a benchmark for our comparison, we use the FOM results, which correspond to a
DNS of the turbulent channel flow.  We compare the three Reg-ROMs, the G-ROM, and
the ROM projection in terms of their accuracy with respect to the FOM benchmark.
We expect the three Reg-ROMs to yield significantly more accurate results than
the standard G-ROM.  We also note that, as mentioned in
Section~\ref{section:rom-filters}, the computational cost of the three Reg-ROMs
is similar, and is on the same order of magnitude as the G-ROM cost.

The rest of this section is organized as follows:
In Section~\ref{section:fom-computational setting}, we present the FOM
computational setting, which is leveraged in
Section~\ref{section:rom-computational setting} to construct the ROMs.
In Section~\ref{section:criteria}, we define the criteria used to evaluate the
ROM performance.
    In Section~\ref{section:offline-online} we outline the efficient
offline-online decomposition of the Reynolds stresses that are used in our
numerical investigation.  We also outline the ROM projection for the Reynolds
stresses, which is used as a benchmark in our numerical investigation.  
Next, we present numerical results of the Reg-ROM comparison for two regimes: In
Section~\ref{section:reproduction-regime}, we present results for the
reproduction regime, i.e., when the ROMs are tested on the training time
interval.  In Section~\ref{section:predictive-regime}, we present results for
the predictive regime, i.e., when the ROMs are tested on a time interval that
is different from the training interval.
Finally, in Section~\ref{section:sensitivity}, we perform a numerical
investigation of the Reg-ROMs' sensitivity with respect to the following
parameters:
(i) the time interval;
(ii) the relaxation parameter, $\chi$; and 
(iii) the filter radius, $\delta$.
The objective of this sensitivity study is to determine which of the three Reg-ROMs 
is more robust with respect to parameter changes.

\subsection{FOM Computational Setting}
    \label{section:fom-computational setting}

\begin{figure}[!ht]
    \centering
        \includegraphics[width=0.55\linewidth]{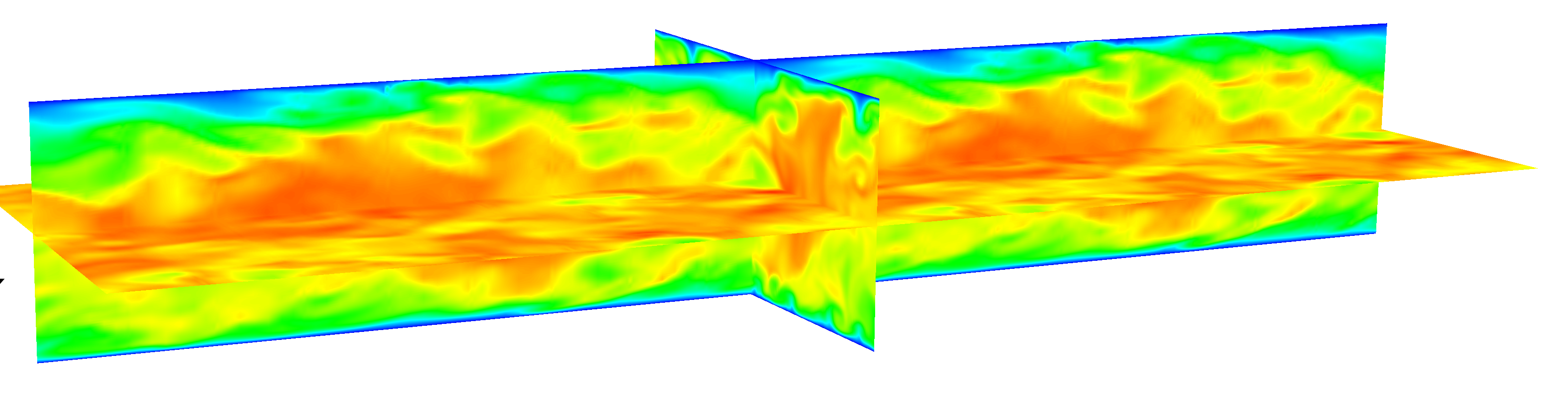}
        \includegraphics[width=0.40\linewidth]{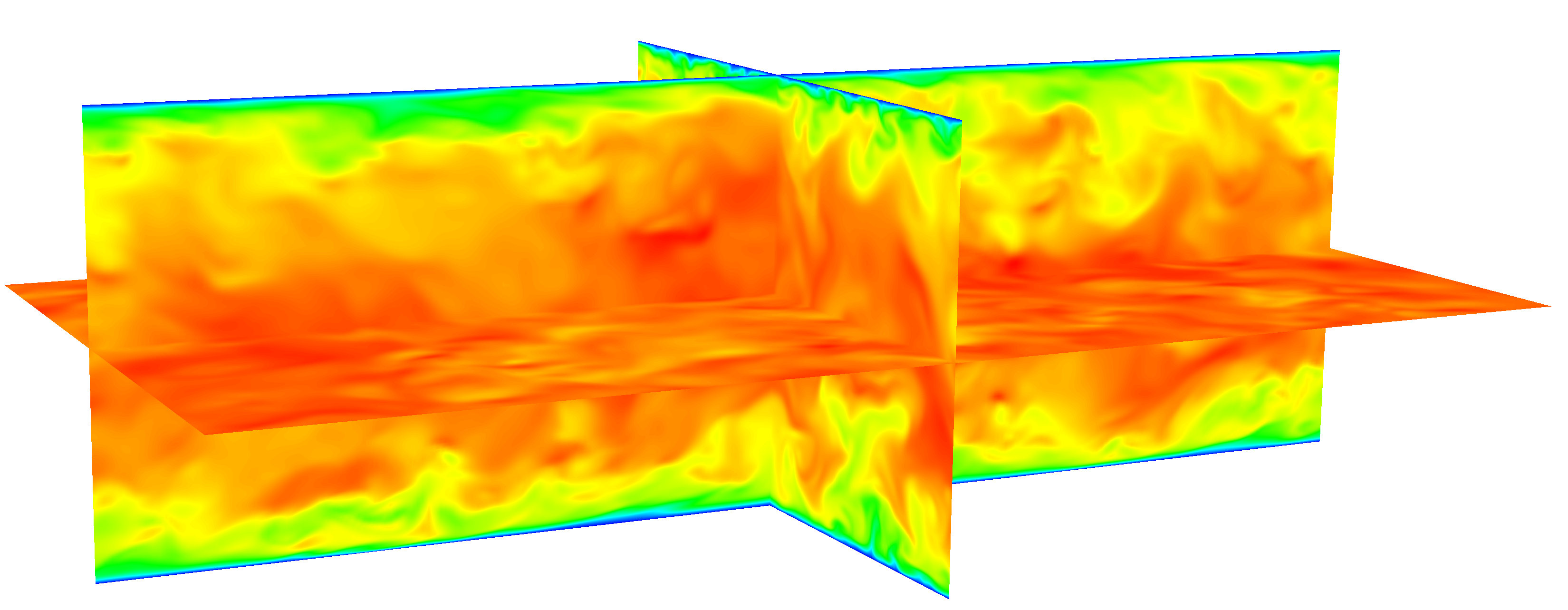}
    \caption{
    Turbulent channel flow: FOM velocity magnitude snapshots 
    at $\rm Re_\tau = 180$ (left) and $395$ (right).} \label{fig:fom_snap}
\end{figure}

In this section, we present the computational setting for the FOM, which has two
main goals: 
(i) to generate the snapshots used in
Section~\ref{section:rom-computational setting} to construct the ROMs; and 
(ii) to serve as a benchmark in the ROM numerical investigation.
Our FOM is a DNS of the turbulent channel flow at
$\rm Re_\tau = 180$ and $\rm Re_\tau = 395$ using the spectral element code Nek5000/RS
\cite{fischer2008nek5000,fischer2022nekrs}. The friction Reynolds number $\rm
Re_\tau$ is based on the friction velocity $u_\tau$ at the wall, channel
half-height $h$, and the fluid kinematic viscosity $\nu$, with $u_\tau =
\sqrt{\tau_w}/\rho$ determined using the wall shear stress, $\tau_w$, and the
fluid density, $\rho$.  FOM velocity magnitude snapshots 
for both $\rm Re_\tau$ are shown in Fig.~\ref{fig:fom_snap}.

For $\rm Re_\tau = 180$, we follow the setup in \cite{kim1987turbulence}, in
which the streamwise (i.e, the $x$-direction) and spanwise (i.e, the
$z$-direction) lengths of the channel are set to $4\pi h$ and $4 \pi h/3$,
respectively, and the channel half-height is set to $h = 1$. 
We consider twice
as many grid points as in \cite{kim1987turbulence}, with $E=5,832$ elements (an
array of $18 \times 18 \times 18$ elements in the $x \times y \times z$
directions), of order $q=9$, for a total of $\cN \approx 4.3$ million grid
points. 
FOM statistics are collected over $1000$ convection time units (CTUs)
and compared against the following two databases: (i) data in
\cite{vreman2014comparison}, which is collected over $3140$ CTUs and has $14.2$
million grid points; and (ii) data in \cite{kim1987turbulence}, which has $2.1$
million grid points \footnote{We couldn't find out how long the statistics are
collected for in \cite{kim1987turbulence}.}.

For $\rm Re_\tau = 395$, we follow the setup in \cite{moser1999direct}, in which
the streamwise and spanwise lengths of the channel are set to $2 \pi$ and $\pi$,
respectively, and the channel half-height is set to $h=1$. We consider twice as
many grid points as in \cite{moser1999direct}, with $E=26,244$ elements (an
array of $36 \times 27 \times 27$ elements in the $x \times y \times z$
directions), of order $q=9$, for a total of $\cN \approx 20$ million grid
points. FOM statistics are collected over $1500$ CTUs and compared against the
data in \cite{moser1999direct}, which has $9.5$ million grid points.

\begin{figure}[h]
    \centering
    \includegraphics[width=0.49\linewidth]{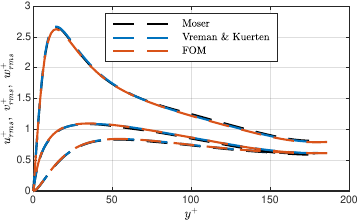}
    \includegraphics[width=0.49\linewidth]{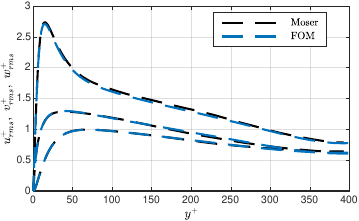}
    \caption{
    FOM 
    2nd order statistics $u^+_{rms},~v^+_{rms}$, and $w^+_{rms}$ validation: 
    (Left) $\rm Re_\tau = 180$, comparison with \cite{kim1987turbulence,vreman2014comparison}. 
    (Right) $\rm Re_\tau = 395$, comparison with 
    \cite{moser1999direct}. 
    }
    \label{fig:fom_stat_compare}
\end{figure}

For both Reynolds numbers, the FOM is run until the solution reaches a
statistically steady state prior to gathering statistics. To validate our FOM,
in Fig.~\ref{fig:fom_stat_compare}, for both $\rm Re_\tau$, we compare the FOM
with the reference data with respect to the 2nd order turbulent statistics
$u^+_{rms}$, $v^+_{rms}$, and $w^+_{rms}$ (in wall-units).  For both $\rm
Re_\tau$, we observe that the FOM 2nd order statistics are in good agreement
with the published results.

\subsection{ROM Computational Setting}
    \label{section:rom-computational setting}

The POD basis functions are constructed using $K=2000$ uniformly distributed snapshots in
the statistically steady state region, which spans $500$ CTUs for $\rm Re_\tau =
180$, and $1000$ CTUs for $\rm Re_\tau = 395$.

For each $\rm Re_\tau$, we compare the ROM performance for different ROM
parameters in both the reproduction and the predictive regimes.  The ROM offline
phase, which includes the construction of the POD basis functions and reduced operators,
is performed using NekROM on the UIUC cluster Delta.  The ROM online phase,
which includes loading the reduced operators and solving the ROM systems, is
performed using Matlab on a workstation.

The implementation of the three Reg-ROMs (\ref{eqn:l-rom}--\ref{eqn:tr-rom}) is
similar: For a given filter order, $m$, and filter radius, $\delta$, all three
Reg-ROMs involve computing the Cholesky factorization of the HOAF
(\ref{eqn:hodf}), that is $(I + \delta^{2m} (B^{-1}A)^m) = R^TR$, and storing the inverse
of the upper triangular matrix, $R^{-1}$. At each time step, a matrix
multiplication $R^{-1}R^{-T}\uu$ is performed to obtain the filtered ROM
coefficients,  $\ubu$. The only difference in the particular Reg-ROM
implementation is how $\ubu$ is used in the Reg-ROM equation. We also emphasize
that, because the most expensive part is to filter the ROM coefficients, the
computational overhead is similar for all three Reg-ROMs.

\subsection{Criteria}
    \label{section:criteria}
    
To evaluate the ROM performance, we use the FOM data as a benchmark, and 
    the streamwise Reynolds normal stress $\langle u'u' \rangle$ and the Reynolds shear
    stress $\langle u'v' \rangle$, which are the two dominant terms in the
    Reynolds stress tensor,
as criteria for accuracy evaluation.
Specifically, we use the following formulas:
\begin{align}
    \varepsilon_{u'u'} \coloneqq \frac{\|\langle u'u' \rangle_{\tr} - \langle u'u' \rangle \|_2}{\|{\langle u'u' \rangle} \|_2},\quad
    \varepsilon_{u'v'} \coloneqq \frac{\|\langle u'v' \rangle_{\tr} - \langle u'v' \rangle \|_2}{\|{\langle u'v' \rangle} \|_2},
    \label{equation:e-reystress}
\end{align}
where $\langle u'u' \rangle$ and $\langle u'v' \rangle$ are the FOM Reynolds
stresses, and $\langle u'u' \rangle_\tr$, $\langle u'v' \rangle_\tr$ are the ROM
Reynolds stresses.

The FOM Reynolds stresses are defined as follows:
\begin{align}
     \langle u'u' \rangle &  \coloneqq \Bigl \langle \left(u -\bracu \right)^2 \Bigr \rangle =
     \Bigl \langle u^2 -2 u \bracu  + \bracu^2 \Bigr \rangle = \Bigl \langle \langle u^2 \rangle - \bracu^2 \Bigr \rangle, \label{eq:u'u'}\\
     \langle u'v' \rangle &  \coloneqq 
     \Bigl \langle  \left( u -\bracu \right) \left( v-\bracv \right) \Bigr \rangle \nonumber \\&=
      \Bigl \langle uv - u \bracv - v \bracu + \bracu \bracv \Bigr \rangle 
     = \Bigl \langle \langle uv \rangle - \bracu \bracv \Bigr \rangle. \label{eq:u'v'}
\end{align}
The ROM Reynolds stresses 
are computed using (\ref{eq:u'u'})-(\ref{eq:u'v'}), but with the ROM approximated solution (\ref{eq:romu}):
\begin{align}
     \langle u'u' \rangle_\tr &  \coloneqq \Bigl \langle \langle u_\tr^2 \rangle - \bracur^2 \Bigr \rangle, \quad
     \langle u'v' \rangle_\tr \coloneqq 
     \Bigl \langle \langle u_\tr v_\tr \rangle - \bracur \bracvr \Bigr \rangle. \label{eq:u'u'_u'v'_direct}
\end{align}
In this paper, an angle bracket $\bigl \langle \cdot \bigr \rangle$ indicates an average over $x$, $z$, and $t$, and is defined as: 
\begin{align}
    \langle u \rangle (y) = \frac{1}{T L_x L_x} \sum_{x,z,t} u(x,y,z,t),
\end{align}
where $T$ is the length of the time interval, $L_x$ is the dimension of the
computational domain in the $x$-direction, $L_z$ is the dimension of the
computational domain in the $z$-direction, $u$ is a scalar field, and a prime
indicates perturbation from this average.

We note that using \cref{eq:u'u'_u'v'_direct} to compute ROM Reynolds stresses
requires accessing the POD basis functions and reconstructing the ROM quantities $\langle
u^2_\tr \rangle$, $\langle u_\tr v_\tr \rangle$, $\bracur^2$, and $\bracur
\bracvr$, which scale with the FOM dimension, $\cN$.  Thus, using
\cref{eq:u'u'_u'v'_direct} to compute $\langle u'u' \rangle_\tr$ and $\langle
u'v' \rangle_\tr$ is inefficient.

\subsubsection{Efficient Offline-Online Reynolds-Stress Evaluation}
    \label{section:offline-online}
To efficiently compute the ROM Reynolds stresses, we use an alternative
approach, based on an offline-online splitting.  
First, we rewrite \cref{eq:u'u'_u'v'_direct} with the POD expansion:
\begin{align}
     \langle u'u' \rangle_\tr & 
     = \Bigl\langle \bigl\langle(\sum^N_{j=0} u_{\tr,j}(t) \varphi_j )(\sum^N_{k=0} u_{\tr,k}(t) \varphi_k) \bigr\rangle -\bigl\langle \sum^N_{j=0}u_{\tr,j}(t) \varphi_j \bigr\rangle^2 \Bigr\rangle, \nonumber \\ &
     = \langle \sum^N_{j=0} \sum^N_{k=0} \langle u_{\tr,j}(t) u_{\tr,k}(t) \varphi_j \varphi_k \rangle - \left( \sum^N_{j=0} \langle u_{\tr,j}(t) \varphi_j \rangle \right) \left( \sum^N_{k=0} \langle u_{\tr,k}(t) \varphi_k \rangle \right) \rangle, \nonumber \\ & = 
      \sum^N_{j=0}\sum^N_{k=0} \langle u_{\tr,j}u_{\tr,k} \rangle_t \langle \varphi_j \varphi_k \rangle_{xz} - \sum^N_{j=0}\sum^N_{k=0} \langle u_{\tr,j} \rangle_t 
      \langle u_{\tr,k} \rangle_t \langle \varphi_j \rangle_{xz} \langle \varphi_k \rangle_{xz} \label{eq:u'u'_off_on},
\end{align}
where $\varphi_j$ and $\varphi_k$ are the POD basis functions in the streamwise
direction, and $\langle \cdot \rangle_t$ and $\langle \cdot \rangle_{xz}$ are
the average operators in $t$ and the $x$-$z$ plane. Then, in the offline stage,
for each $y$ coordinate, we compute and store $\langle \varphi_j \varphi_k
\rangle_{xz}$ and $\langle \varphi_j \rangle_{xz}$ for all $j,~k = 0,\ldots,N$. 
Finally, in the online stage, we compute $\langle {u_{\tr,j} u_{\tr,k}}
\rangle_t$ and $\langle {u_{\tr,j}} \rangle_t$ for all $j,~k=0,\ldots,N$.  Thus,
to construct the ROM streamwise Reynolds normal stress $\langle u'u' \rangle_\tr$ at a
given point $y^*$, we use (\ref{eq:u'u'_off_on}), which is independent of $\cN$
and, thus, does not significantly increase the computational cost. The
offline-online splitting for the ROM Reynolds shear stress $\langle u'v'
\rangle_\tr$ can be derived similarly.

We also assess the ROM performance with the ROM-projection Reynolds stresses
$\langle u'u' \rangle_{\tProj}$ and $\langle u'v' \rangle_{\tProj}$:
\begin{align}
     \langle u'u' \rangle_\tProj &  \coloneqq \Bigl \langle \langle u_\tProj^2 \rangle - \bracuproj^2 \Bigr \rangle, \label{eq:u'u'_Proj}\\
     \langle u'v' \rangle_\tProj & \coloneqq 
     \Bigl \langle \langle u_\tProj v_\tProj \rangle - \bracuproj \bracvproj \Bigr \rangle, \label{eq:u'v'_POD}
\end{align}
where the 
projected solution onto the reduced space is defined as:
\begin{equation} \bu_\tProj = \bphi_0 + \sum^N_{j=0} \tilde{u}_{\tProj,j}
\bphi_j(\bx), ~\text{and}~\tilde{u}_{\tProj,j} = (\bphi_j, \bu - \bphi_0).
\end{equation}
    The ROM-projection Reynolds stress can be computed using
    (\ref{eq:u'u'_off_on}) but with $\langle \tilde{u}_{\tProj,j}
    \tilde{u}_{\tProj,j} \rangle_t$ and $\langle \tilde{u}_{\tProj,j} \rangle_t$
    quantities, and is used to compute the error $\euv$
    (\ref{equation:e-reystress}). The ROM projection represents the best
    theoretical approximation of the training data in the given ROM space, and we
    use it as a benchmark to assess the three Reg-ROMs.
\subsection{Reproduction Regime} 
    \label{section:reproduction-regime}

In this section, we perform a numerical investigation of the three Reg-ROMs:
L-ROM (Section~\ref{section:l-rom}), EFR-ROM
(Section~\ref{section:efr-rom}), and the new TR-ROM
(Section~\ref{section:tr-rom}) 
    in the reproduction regime, i.e., in the time interval in which the
    snapshots were collected, at $\rm Re_\tau=180$
    (Section~\ref{section:reproduction-regime-180}) and $\rm Re_\tau=395$
    (Section~\ref{section:reproduction-regime-395}).
For comparison purposes, we include results for the G-ROM
(Section~\ref{section:g-rom}) and the ROM projection, which represents the best
theoretical approximation of the FOM solution in the given ROM space.  The
Reg-ROM accuracy is expected to be significantly higher than the G-ROM accuracy.
To quantify the ROM accuracy, we use $\euu$ and $\euv$,  which are the relative
$\ell^2$ errors of the streamwise Reynolds normal stress $\langle u'u' \rangle $
and Reynolds shear stress $\langle u'v' \rangle$ defined in
\cref{equation:e-reystress}.

In our numerical investigation, 
we also consider the following parameters:
    For all the ROMs, we utilize $10$ values for the ROM dimension, $N \in \{10,20,\ldots,$
    $90,100\}$. 
    For each Reg-ROM, 
        we use four HOAF orders: $m=1$ (which corresponds to the classical
        DF~\eqref{equation:df-linear-system}) and $m = 2, 3, 4$ (which
        correspond to the HOAF~\eqref{eqn:hodf}).
        %
        For each $N$ and $m$ value, the values of the filter radius, $\delta$,
        are chosen as follows: 
        For EFR-ROM and TR-ROM, $10$, $25$, and $10$
        $\delta$ values are uniformly sampled from the intervals
        $[0.001,~0.01]$, $[0.01, 0.1]$, and $[0.1, 1]$, respectively.  This
        yields a total of $43$ values for $\delta$ from the interval
        $[0.001,~1]$.
        In addition, we choose $4$ uniformly sampled values for the
        relaxation parameter, $\chi$, in the interval
        $[\Delta t = 0.005,~1]$.  We choose this $\chi$ range because $\chi =
        \cO(\Delta t)$ is commonly used in EFR-ROM
        simulations~\cite{ervin2012numerical,strazzullo2022consistency}. 
        For L-ROM, we uniformly sample $15$ additional $\delta$ values from the
        interval $[0.1,~0.2]$, which yields a total of $56$ values for $\delta$
        from the interval $[0.001,~1]$ for $\rm Re_\tau=180$. 
        For $\rm Re_\tau = 395$, $30$ additional $\delta$ values are uniformly
        sampled from the interval $[0.1,~0.3]$, which yields a total of $69$
        values for $\delta$.  

We emphasize that, in our Reg-ROM numerical investigation, we use four
parameters: the ROM dimension, $N$, the filter order, $m$, the filter radius,
$\delta$, and (for EFR-ROM and TR-ROM) the relaxation parameter, $\chi$.  Thus,
to ensure a clear comparison of the three Reg-ROMs, we fix the $\delta$ and
$\chi$ parameters to their optimal values (i.e., the values that yield the most
accurate Reynolds shear stress $\langle u'v' \rangle$ for each Reg-ROM), and
plot $\euv$ for all the parameter values for $N$ and $m$.  We note that we show
only $\euv$ results for the following two reasons: (i) We find that $\euu$ behaves
similarly to $\euv$.  Specifically, $\euu$ is generally smaller than $\euv$ at
$\rm Re_\tau=180$, and similar to $\euv$ at
$\rm Re_\tau=395$. (ii) In turbulent channel flow investigations, the
approximation of the Reynolds shear stress $\langle u'v' \rangle$ is more
challenging than the approximation of streamwise Reynolds normal stress $\langle
u'u' \rangle$.  For the $\euu$ results, we refer reader to the dissertation
\cite{tsai2023pmordevelopment}.  For completeness, in the sensitivity study in
Section~\ref{section:sensitivity}, we include Reg-ROM results for all the
$\delta$ and $\chi$ values. 

Finally, in Section~\ref{section:reproduction-regime-summary}, we present a
summary of the Reg-ROM comparison in the reproduction regime.

\subsubsection{$\rm Re_\tau= 180$}
    \label{section:reproduction-regime-180}

\begin{figure}[!ht]
    \centering
    \begin{subfigure}{0.49\columnwidth}
    \caption{G-ROM}
    \label{subfig:recon_180_grom}
    \includegraphics[width=1\textwidth]{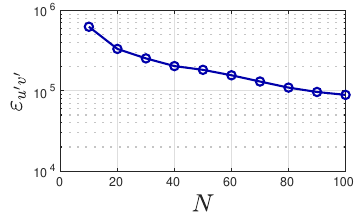}
    \end{subfigure}
    \begin{subfigure}{0.49\columnwidth}
    \caption{L-ROM}
    \label{subfig:recon_180_lrom}
    \includegraphics[width=1\textwidth]{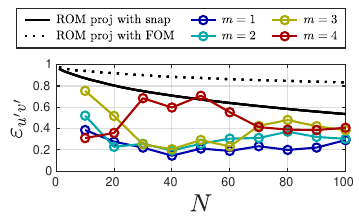}
    \end{subfigure}\\
    \begin{subfigure}{0.49\columnwidth}
    \caption{EFR-ROM}
    \label{subfig:recon_180_efr}
    \includegraphics[width=1\textwidth]{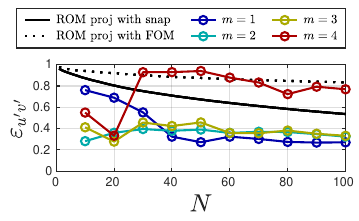}
    \end{subfigure}
    \begin{subfigure}{0.49\columnwidth}
    \caption{TR-ROM}
    \label{subfig:recon_180_tr}
    \includegraphics[width=1\textwidth]{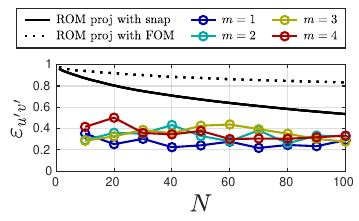}
    \end{subfigure}
    \caption{ The relative error $\euv$~\cref{equation:e-reystress} of G-ROM,
    ROM projection, L-ROM, EFR-ROM, and TR-ROM in the \textit{reproduction
    regime} at $\rm Re_\tau =180$ for different $N$ and $m$ values with optimal
    $\delta$ and $\chi$ values.}
    \label{fig:Re180_recon_conv}
\end{figure}

In \cref{fig:Re180_recon_conv}, we plot the relative $\ell^2$ error $\euv$
\cref{equation:e-reystress} for different ROM dimensions, $N$, and filter
orders, $m$, for the G-ROM, the ROM projection, and the three Reg-ROMs at $\rm
Re_\tau= 180$. 
Two sets of ROM-projection results are shown: The solid and the
dotted lines represent 
    the errors of $\langle u'v' \rangle_\tProj$
    (\ref{eq:u'v'_POD}) for 
    the snapshot data and the FOM data, respectively. 
It is expected that the error $\euv$ for the snapshot data is smaller than the
error for the FOM data, because the POD basis set is constructed using only
$K=2000$ snapshots instead of using all the FOM solutions.  For each Reg-ROM,
the error is plotted for the optimal $\delta$ values and, for EFR-ROM and
TR-ROM, for the optimal $\chi$ values.

\cref{subfig:recon_180_grom} displays the G-ROM results.  This plot shows that,
for all $N$ values, the G-ROM results are very inaccurate.  Even with $N=100$,
G-ROM fails to reconstruct $\langle u'v' \rangle$ with an error of $\cO(10^5)$.

\cref{subfig:recon_180_lrom} displays the L-ROM results for each $N$ and $m$ with the optimal $\delta$, 
along with the results of the ROM projection for comparison purposes.
For $N \ge 30$, $m=1$ yields the most accurate results, achieving an error of $15\%$ for $N=40$.
Conversely, for $N \le 20$, a higher-order filter yields better results. Specifically, $m=4$ achieves an error of $31\%$ for $N=10$ and $m=2$ achieves an error of $23\%$ for $N=20$.
With the exception of $(N, m) = (50, 4)$, L-ROM is more accurate than the ROM projection.

\cref{subfig:recon_180_efr} displays the EFR-ROM results for each $N$ and $m$ with the optimal $\chi$ and $\delta$ values, along with the results of the ROM projection for comparison purposes. 
For $N \ge 40$, $m=1$ yields the most accurate results, achieving an error of $27\%$ for $N=90$.
Conversely, for $N \le 30$, a higher-order filter yields better results. Specifically,  $m=2$ achieves an error of $28\%$ and $40\%$ for $N=10$ and $N=30$, respectively. Additionally, $m=3$ achieves an error of  $45\%$ for $N=20$.
EFR-ROM is more accurate than the ROM projection for all $N$ with $m=1,~2,~3$. With $m=4$, it is only more accurate than the ROM projection for $N\le 30$.

\cref{subfig:recon_180_tr} displays the TR-ROM results for each $N$ and $m$ with the optimal $\chi$ and $\delta$ values, along with the results of the ROM projection for comparison purposes.
For $20 \le N \le 90$, $m=1$ yields the most accurate results, achieving an error of $22\%$ for $N=70$.
For $N = 10$ and $N=100$, $m=3$ yields the most accurate results, achieving an error of $29\%$ and $28\%$, respectively. In addition, TR-ROM is more accurate than the ROM projection for all $N$ and $m$ values.

\subsubsection{$\rm Re_\tau= 395$}
    \label{section:reproduction-regime-395}

\begin{figure}[!ht]
    \centering
    \begin{subfigure}{0.49\columnwidth}
    \caption{G-ROM}
    \label{subfig:recon_395_grom}
    \includegraphics[width=1\textwidth]{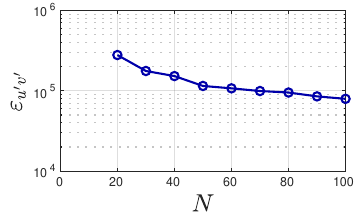}
    \end{subfigure}
    \begin{subfigure}{0.49\columnwidth}
    \caption{L-ROM}
    \label{subfig:recon_395_lrom}
    \includegraphics[width=1\textwidth]{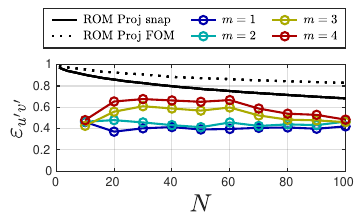}
    \end{subfigure}\\
    \begin{subfigure}{0.49\columnwidth}
    \caption{EFR-ROM}
    \label{subfig:recon_395_efr}
    \includegraphics[width=1\textwidth]{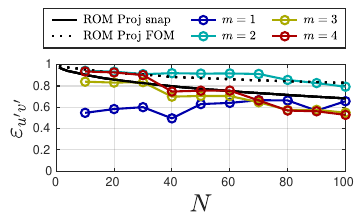}
    \end{subfigure}
    \begin{subfigure}{0.49\columnwidth}
    \caption{TR-ROM}
    \label{subfig:recon_395_tr}\includegraphics[width=1\textwidth]{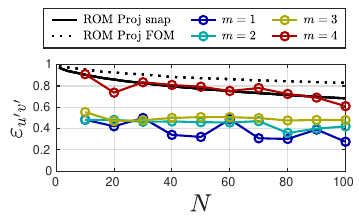}
    \end{subfigure}
    \caption{ The relative error $\euv$~\cref{equation:e-reystress} of G-ROM,
    ROM projection, L-ROM, EFR-ROM, and TR-ROM in the \textit{reproduction
    regime} at $\rm Re_\tau = 395$ for different $N$ and $m$ values with optimal
    $\delta$ and $\chi$ values.}
    \label{fig:Re395_recon_conv}
\end{figure}

In \cref{fig:Re395_recon_conv}, we plot the relative $\ell^2$ error $\euv$
\cref{equation:e-reystress} for different $N$ and $m$ values for the G-ROM, the
ROM projection, and the three Reg-ROMs at $\rm Re_\tau= 395$. 
Two sets of ROM-projection results 
    for the snapshot data and the FOM data are shown.
Note that the error $\euv$ for the snapshot data at $\rm
Re_\tau=395$ is larger (approximately $70\%$) 
than the corresponding results
at $\rm Re_\tau=180$ (approximately $52\%$). This difference is expected because
the $\rm Re_\tau=395$ solution 
is more turbulent and, as a result,
requires a 
larger number of modes to achieve a satisfactory approximation.
For each Reg-ROM, the error is plotted for the optimal $\delta$
values and, for EFR-ROM and TR-ROM, for the optimal $\chi$ values.

\cref{subfig:recon_395_grom} displays the G-ROM results.  Just as in
Section~\ref{section:reproduction-regime-180}, for all $N$ values, the G-ROM
results are very inaccurate. Even with $N=100$, $\euv$ is still about
$\cO(10^5)$.

\cref{subfig:recon_395_lrom} displays the L-ROM results for each $N$ and $m$
with the optimal $\delta$ values, along with the results of the ROM projection
for comparison purposes.  For all $N$ and $m$ values, the error is 
much higher 
than the error for $\rm Re_\tau = 180$.  For $N\geq 20$,
the most accurate results are achieved with $m=1$, with an error of $37\%$ for
$N=40$. For $N=10$, $m=3$ achieves an error of $43\%$.  
Compared to the ROM projection, L-ROM is more accurate for all
values of $N$ and $m$. 

\cref{subfig:recon_395_efr} displays the EFR-ROM results for each $N$ and $m$
with the optimal $\chi$ and $\delta$ values, along with the results of the ROM
projection for comparison purposes.  These results are qualitatively different
from the EFR-ROM results for $\rm Re_\tau = 180$.  For $N \le 60$, $m=1$ yields
the most accurate results and achieves an error of $50\%$ for $N=40$. For $N \ge
70$, higher-order filter yields better results. Specifically, $m=4$ achieves an
error of $53\%$ for $N=100$.  Moreover, EFR-ROM is found to be more accurate
than the ROM projection for all $N$ with $m=1$ and for $N\ge 60$ with $m=3,~4$.
With $m=2$, EFR-ROM 
has a similar level of accuracy as the ROM projection.

\cref{subfig:recon_395_tr} displays the TR-ROM results for each $N$ and $m$ with
the optimal $\chi$ and $\delta$ values, along with the results of the ROM
projection for comparison purposes.  For almost all $N$ values, $m=1$ yields the
most accurate results, achieving an error of around $28\%$ for $N=100$.  These
results also show that $m=2$ and $m=3$ yield similar results, while $m=4$ is
found to be the least accurate.  Moreover, TR-ROM is found to be more accurate
than the ROM projection for all $N$ and $m$ values.

\subsubsection{Summary}
    \label{section:reproduction-regime-summary}

Overall, our numerical investigation in the reproduction regime yields the following general conclusions:

All three Reg-ROMs are significantly more accurate than the standard G-ROM. 
    In fact, with respect to several second-order turbulence statistics, the
    errors of the three Reg-ROMs equipped with carefully tuned spatial filtering
    are much lower than the projection error.

Finally, our numerical investigation demonstrates that, for $\rm Re_\tau =180$,
all three Reg-ROMs with $m=1$ (i.e., low-order filtering) consistently produce
the most accurate results for large $N$ values, while a higher-order filter is
more effective for low $N$ values. For $\rm Re_\tau =395$, L-ROM and TR-ROM with
$m=1$ yield the most accurate results for all $N$ values, 
    while EFR-ROM yields the most accurate results for low $N$ values with $m=1$, 
    and for high $N$ values with $m=3,~4$.  

To facilitate the comparison of the three Reg-ROMs, in
\cref{table:accuracy-ranking}, we rank them based on the error $\euv$ and $\euu$
achieved for the $N$ and $m$ values investigated.  Specifically, for both
Reynolds numbers, we list the Reg-ROMs' rank, the lowest $\euv$, the
corresponding $\euu$, the ROM dimension $N$ and the filter order $m$ for which
the lowest $\euv$ is achieved.  The results in \cref{table:accuracy-ranking}
yield the following conclusions:

For $\rm Re_\tau=180$, TR-ROM is the most accurate model with $\euv \approx
22\%$ and $\euu \approx 16\%$, followed by EFR-ROM and L-ROM.  For $\rm
Re_\tau=395$, TR-ROM is still the most accurate model with $\euv \approx 28\%$
and $\euu \approx 26\%$, followed by L-ROM and EFR-ROM.  In addition, we find
$\euu$ is smaller compared to $\euv$ except for L-ROM for $\rm Re_\tau =180$ and
$\euu$ has similar level of accuracy as $\euv$ for $\rm Re_\tau =395$.  Moreover,
the results in \cref{table:accuracy-ranking} also show that $m=1$ (i.e.,
low-order filtering) yields the most accurate results.  Finally, these results
show that TR-ROM requires large $N$ values to achieve its best accuracy, whereas
L-ROM yields best accuracy with small $N$. For EFR-ROM, large $N$ is required
for $\rm Re_\tau=180$, and small $N$ is required for $\rm Re_\tau =395$.

\begin{table}[!ht]
\caption{Reg-ROM accuracy ranking in the reproduction regime for $\rm Re_\tau = 180$ (top
rows) and $395$ (bottom rows).  The following parameters are listed: Reg-ROMs'
rank, the lowest $\euv$, the corresponding $\euu$, the ROM dimension $N$ and the
filter order $m$ for which the lowest $\euv$ is achieved.
\label{table:accuracy-ranking}}
\centering
\begin{tabular}{|l|l|l|l|}
\hline
$\rm Re_\tau=180$  & L-ROM & EFR-ROM & TR-ROM \\ \hline
Rank & 3 & 2       & 1      \\ \hline
$\euv$ & $\approx 19\%$    & $\approx 27\%$       & $\approx 22\%$      \\ \hline
$\euu$ & $\approx 37.5\%$    & $\approx 13\%$       & $\approx 16\%$      \\ \hline
N & $40$    & $90$ & $70$      \\ \hline
$m$ &  1    & 1       & 1     \\ \hline \hline
$\rm Re_\tau=395$ & L-ROM & EFR-ROM & TR-ROM \\ \hline
Rank & 2     & 3       & 1      \\ \hline
$\euv$ & $\approx 37\%$    & $\approx 50\%$       & $\approx 28\%$      \\ \hline
$\euu$ & $\approx 40\%$    & $\approx 40\%$       & $\approx 26\%$      \\ \hline
N & $20$    & $40$ & $100$      \\ \hline
$m$ & 1    & 1       & 1     \\ \hline
\end{tabular}
\end{table}

\begin{figure}[!ht]
    \centering
    \begin{subfigure}{0.49\columnwidth}
    \caption{L-ROM, $\rm Re_\tau =180$}
    \label{subfig:}
    \includegraphics[width=1\textwidth]{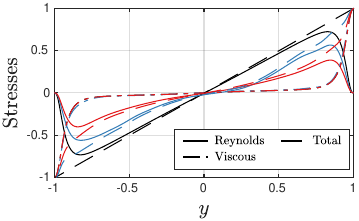}
    \end{subfigure}
    \begin{subfigure}{0.49\columnwidth}
    \caption{L-ROM, $\rm Re_\tau=395$}
    \label{subfig:}
    \includegraphics[width=1\textwidth]{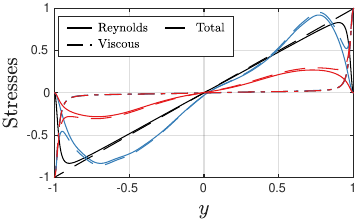}
    \end{subfigure}\\
    \begin{subfigure}{0.49\columnwidth}
    \caption{EFR-ROM, $\rm Re_\tau=180$}
    \label{subfig:}
    \includegraphics[width=1\textwidth]{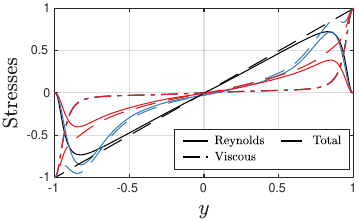}
    \end{subfigure}
    \begin{subfigure}{0.49\columnwidth}
    \caption{EFR-ROM, $\rm Re_\tau =395$}
    \label{subfig:}
    \includegraphics[width=1\textwidth]{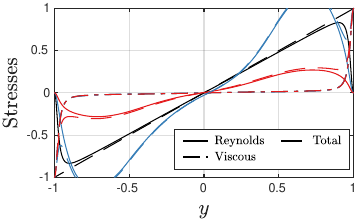}
    \end{subfigure}\\
    \begin{subfigure}{0.49\columnwidth}
    \caption{TR-ROM, $\rm Re_\tau=180$}
    \label{subfig:}
    \includegraphics[width=1\textwidth]{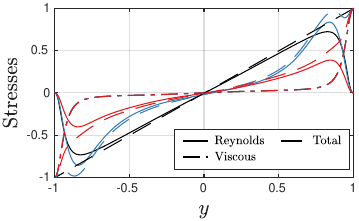}
    \end{subfigure}
    \begin{subfigure}{0.49\columnwidth}
    \caption{TR-ROM, $\rm Re_\tau=395$}
    \label{subfig:}
    \includegraphics[width=1\textwidth]{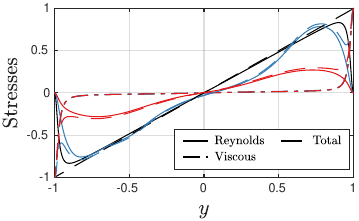}
    \end{subfigure}
    \caption{Comparison of the total, viscous, and Reynolds shear stresses of
    the FOM (black), the ROM projection (red), and the optimal Reg-ROMs, 
    listed in Table~\ref{table:accuracy-ranking} (blue) in the \textit{reproduction
    regime} for $\rm Re_\tau =180$ (left) and $\rm Re_\tau=395$ (right). The
    total shear stress is the sum of the viscous and the Reynolds shear stress.}
    \label{fig:recon_profiles}
\end{figure}

In Fig.~\ref{fig:recon_profiles}, we compare the total, viscous, and Reynolds
shear stresses of the optimal Reg-ROMs (listed in
Table~\ref{table:accuracy-ranking}) along with the results of the FOM and the
ROM projection in the reproduction regime for $\rm Re_\tau=180$ and $\rm
Re_\tau=395$.  The total shear stress $\tau(y)$ is the sum of the viscous shear
stress $\rho \, \nu \, d\langle u \rangle/dy$ and the Reynolds shear stress
$-\rho \, \langle uv \rangle$, and its distribution is linear
\cite{pope2000turbulent}.  In each model, the three shear stresses are
normalized with the model's wall shear stress $\tau_w \equiv \rho \, \nu\,
\left(d\langle u \rangle/dy\right)_{y=-1}$.

In terms of the viscous shear stress, the results of the three Reg-ROMs are in
good agreement with those of the FOM and the ROM projection for both $\rm
Re_\tau$.
In terms of the Reynolds shear stress, for $\rm Re_\tau =180$, we observe that
the L-ROM result 
is smaller than the FOM result. 
On the other hand, 
the results of the EFR-ROM and the TR-ROM are similar, and both have higher values
in the boundary layer and lower values outside the boundary layer compared to
the FOM.  Notice that both the EFR-ROM and the TR-ROM are able to capture the
slope of the Reynolds shear stress. For $\rm Re_\tau =395$, we find 
that the result of the TR-ROM is the best, followed by L-ROM and EFR-ROM.
Moreover, although the results of the Reg-ROMs are not perfect, we find that the
results are much better than the ROM projection. This also indicates that 
$N=100$ POD bases are not able to reconstruct the Reynolds stress.
Finally, as a result of the discrepancy in the Reynolds shear stress, the total
shear stress is not linear in all three Reg-ROMs for both $\rm Re_\tau$.

\subsection{Predictive Regime} 
    \label{section:predictive-regime}

In this section, we perform a numerical investigation of the three Reg-ROMs:
L-ROM (Section~\ref{section:l-rom}), EFR-ROM (Section~\ref{section:efr-rom}),
and the new TR-ROM (Section~\ref{section:tr-rom}) for the \textit{predictive
regime},
    on a time interval that is $500$ CTUs larger than the 
    time
    interval on which snapshots were collected, at $\rm Re_\tau= 180$
    (Section~\ref{section:predictive-regime-180}) and $\rm
    Re_\tau= 395$ (Section~\ref{section:predictive-regime-395}).  Hence, $1000$
    CTU and $1500$ CTU time windows are considered for $\rm
    Re_\tau = 180$ and $395$, respectively. 
For comparison purposes, we include results for the G-ROM
(Section~\ref{section:g-rom}) and the ROM projection. 
To quantify the ROM accuracy, we use $\euu$ and $\euv$, which are the relative
$\ell^2$ error of $\langle u'u' \rangle$ and $\langle u'v' \rangle$ defined in \cref{equation:e-reystress}. 
    In order to measure the accuracy of the predicted $\langle u'u' \rangle $ and $\langle u'v' \rangle$, an
    additional $500$ CTUs of FOM simulations are performed for both $\rm
    Re_\tau$.
For the $\euu$ results, we refer reader to the dissertation \cite{tsai2023pmordevelopment}.

In our numerical investigation, we use the same parameter values for $N$ and $m$
as the values used in Section~\ref{section:reproduction-regime}.  
    For each $N$ and $m$ values, we plot $\euv$ with $(\delta,~\chi)_{\trecon}$,
    that is, the $\delta$ and $\chi$ parameter values that were found to be
    optimal in the reproduction regime in
    Section~\ref{section:reproduction-regime}.
We emphasize that our strategy is different from that used in
Section~\ref{section:reproduction-regime}: Instead of optimizing the $\delta$
and $\chi$ values on the entire predictive time interval, we use the values that
were optimized over the shorter time interval of the reproduction regime.  Thus,
in this section, we are investigating the predictive power of both the Reg-ROMs
{\it and} their associated parameters.

Finally, in Section~\ref{section:predictive-regime-summary}, we present a
summary of the Reg-ROM comparison in the predictive regime.

\subsubsection{$\rm Re_\tau= 180$}
    \label{section:predictive-regime-180}
    
\begin{figure}[!ht]
    \centering
    \begin{subfigure}{0.49\columnwidth}
    \caption{G-ROM}
    \label{subfig:pred_180_grom}
    \includegraphics[width=1\textwidth]{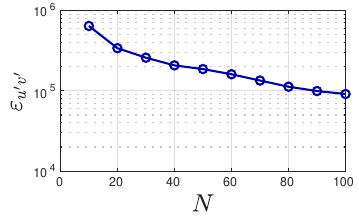}
    \end{subfigure}
    \begin{subfigure}{0.49\columnwidth}
    \caption{L-ROM}
    \label{subfig:pred_180_lrom}
    \includegraphics[width=1\textwidth]{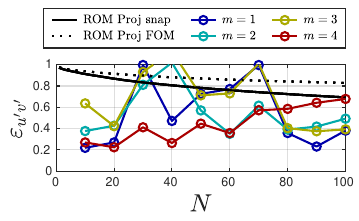}
    \end{subfigure}\\
    \begin{subfigure}{0.49\columnwidth}
    \caption{EFR-ROM}
    \label{subfig:pred_180_efr}
    \includegraphics[width=1\textwidth]{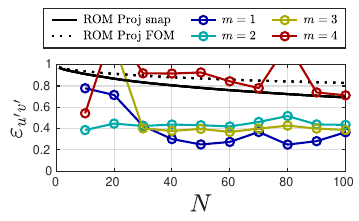}
    \end{subfigure}
    \begin{subfigure}{0.49\columnwidth}
    \caption{TR-ROM}
    \label{subfig:pred_180_tr}
    \includegraphics[width=1\textwidth]{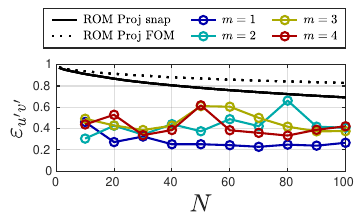}
    \end{subfigure}
    \caption{
    The relative error $\euv$~\cref{equation:e-reystress} of G-ROM, ROM projection, L-ROM, EFR-ROM, and TR-ROM in the \textit{predictive regime} at $\rm Re_\tau =180$ for different $N$ and $m$ values with optimal $\delta$ and $\chi$ values \textit{optimized in the reproduction regime}. 
    }
    \label{fig:Re180_pred_conv}
\end{figure}

In \cref{fig:Re180_pred_conv}, we plot the relative $\ell^2$ error $\euv$
\cref{equation:e-reystress} for different $N$ and $m$ values for the G-ROM, the
ROM projection, and the three Reg-ROMs at $\rm Re_\tau= 180$.  We emphasize
that, to test the predictive capabilities of the Reg-ROM parameters, we plot the
error for $(\delta,~\chi)_\trecon$ that were optimized in the
\textit{reproduction regime} (Section \ref{section:reproduction-regime-180}).

\cref{subfig:pred_180_grom} displays the G-ROM results. As in the reproduction
regime, for all $N$ values, the G-ROM results are very inaccurate. 

\cref{subfig:pred_180_lrom} displays the L-ROM results for each $N$ and $m$ with
$\delta_\trecon$ along with the ROM projection results for comparison purposes. 
For $N = 10$ and $N\ge 80$, $m=1$ yields the most accurate results, achieving an
error of $22\%$ for $N=10$. For $20 \le N \le 70$, $m=4$ yields the most
accurate results, achieving an error of $22\%$ for $N=20$.
For the majority of $N$ and $m$ values, except for the $(m, N)$ values 
$(1, 30)$, $(1, 70)$, $(2, 40)$, $(3, 30)$, $(3, 40)$, and $(3, 70)$, L-ROM is
more accurate than the ROM projection. 

However, $\euv$ exhibits a higher sensitivity to changes in $N$ across all $m$
values compared to the results observed in the reproduction regime (as discussed
in Section \ref{section:reproduction-regime-180}). 

\cref{subfig:pred_180_efr} displays the EFR-ROM results for each $N$ and $m$
with $(\delta,~\chi)_\trecon$ along with the ROM projection results for
comparison purposes.  For $N \le 30$, higher-order filter yields better results.
    Specifically, for $N=10$ and $N=20$, EFR-ROM achieves an error of $39\%$ and
    $45\%$, respectively, with $m=2$.  Additionally, for $N=30$, EFR-ROM
    achieves an error of $40\%$ with $m=3$.
For $N \ge 40$, $m=1$ yields the most accurate results, achieving an error of
$25\%$ for $N=80$.
Compared to the ROM projection, EFR-ROM is more accurate for $m\le 3$,
except for $N=10$ and $N=20$ with $m=3$. For $m=4$, EFR-ROM is not accurate and
its level of accuracy is similar to that of the ROM projection.

\cref{subfig:pred_180_tr} displays the TR-ROM results for each $N$ and $m$ with
$(\delta,~\chi)_\trecon$ along with the ROM projection results for comparison
purposes. 
For $N=10$, $m=2$ yields the lowest error of $30\%$.  For $N \ge 20$, $m=1$
yields the best results, achieving an error of $23\%$ for $N=70$.  TR-ROM is
more accurate than the ROM projection for all $N$ and $m$ values.

\subsubsection{$\rm Re_\tau= 395$}
    \label{section:predictive-regime-395}

\begin{figure}[!ht]
    \centering
    \begin{subfigure}{0.49\columnwidth}
    \caption{G-ROM}
    \label{subfig:pred_395_grom}
    \includegraphics[width=1\textwidth]{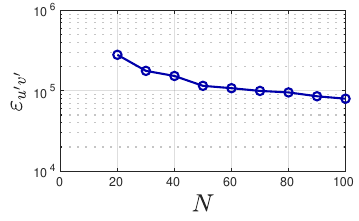}
    \end{subfigure}
    \begin{subfigure}{0.49\columnwidth}
    \caption{L-ROM}
    \label{subfig:pred_395_lrom}
    \includegraphics[width=1\textwidth]{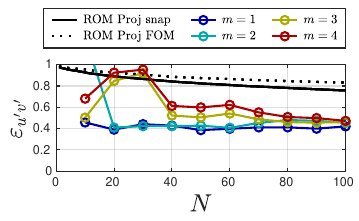}
    \end{subfigure}\\
    \begin{subfigure}{0.49\columnwidth}
    \caption{EFR-ROM}
    \label{subfig:pred_395_efr}
    \includegraphics[width=1\textwidth]{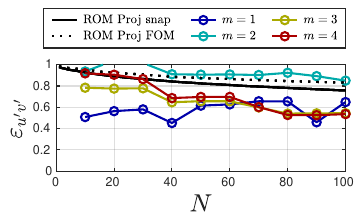}
    \end{subfigure}
    \begin{subfigure}{0.49\columnwidth}
    \caption{TR-ROM}
    \label{subfig:pred_395_tr}
    \includegraphics[width=1\textwidth]{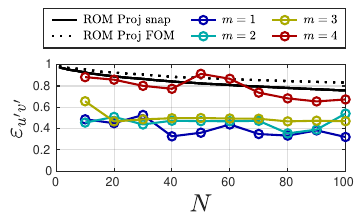}
    \end{subfigure}
    \caption{ The relative error $\euv$~\cref{equation:e-reystress} of G-ROM,
    ROM projection, L-ROM, EFR-ROM, and TR-ROM in the \textit{predictive regime}
    at $\rm Re_\tau =395$ for different $N$ and $m$ values with optimal $\delta$
    and $\chi$ values \textit{optimized in the reproduction regime}. 
    }
    \label{fig:Re395_pred_conv}
\end{figure}

In \cref{fig:Re395_pred_conv}, we plot the relative $\ell^2$ error $\euv$ for
different $N$ and $m$ values for the G-ROM, the ROM projection, and the three
Reg-ROMs at $\rm Re_\tau= 395$.  We emphasize that, to test the predictive
capabilities of the Reg-ROM parameters, we plot the error for
$(\delta,~\chi)_\trecon$ that were optimized in the \textit{reproduction regime}
(Section \ref{section:reproduction-regime-395}).

\cref{subfig:pred_395_grom} displays the G-ROM results.  Just as in
Section~\ref{section:predictive-regime-180}, for all $N$ values, the G-ROM
results are very inaccurate with an error of $\cO(10^5)$ for $N=100$.

\cref{subfig:pred_395_lrom} displays the L-ROM results for each $N$ and $m$ with
$\delta_\trecon$, along with the ROM projection results for comparison purposes.
For almost all $N$ values, except for $N=30$ and $N=40$, $m=1$ yields the most
accurate results, achieving an error of $39\%$ with $N=20$.  For $m=3$ and
$m=4$, the accuracy of the L-ROM is improved as $N$ increases but is still
larger than the $m=1$ case.  
Compared to the ROM projection, L-ROM is more
accurate except for the $(m,N)$ pairs $(2,10), (3,30)$, $(4,20)$, and $(4, 30)$.

\cref{subfig:pred_395_efr} displays the EFR-ROM results for each $N$ and $m$
with $(\delta,~\chi)_\trecon$ along with the ROM projection results for
comparison purposes.  For $N \le 60$, $m=1$ yields the most accurate results,
achieving an error of $45\%$ with $N=40$.  For $N \ge 70$, except for $N=90$, a
higher-order filter yields better results. Specifically, $m=4$ achieves an error
of $53\%$ with $N=80$. For $N=90$, $m=1$ achieves an error of $46\%$.
Compared to the ROM projection, EFR-ROM is
more accurate for all $N$ values and $m=1,~3,~4$. For $m=2$, the EFR-ROM's 
level of accuracy is similar to that of the ROM projection.

\cref{subfig:pred_395_tr} displays the TR-ROM results for each $N$ and $m$ with
$(\delta,~\chi)_\trecon$, along with the ROM projection results for comparison
purposes.  For all $N$ values except $N=10$ and $N=30$, $m=1$ yields the most
accurate results, achieving an error of around $32\%$ for $N=100$.  These
results also show that $m=2$ and $m=3$ yield similar accuracy, and $m=4$ is the
least accurate. TR-ROM is more accurate than the ROM projection for all $N$ and
$m$ values, except for the $(m, N)$ pairs $(4, 50)$ and $(4, 60)$.

\subsubsection{Summary}
    \label{section:predictive-regime-summary}
    
Overall, our numerical investigation in the predictive regime yields the following general conclusions:

All three Reg-ROMs are significantly more accurate than the standard G-ROM. 

In fact, with respect to several second-order turbulence statistics, 
the errors of the three Reg-ROMs equipped with carefully tuned spatial filtering 
are much lower than 
the projection error. 

Finally, our numerical investigation demonstrates that, for $\rm Re_\tau =180$,
EFR-ROM and TR-ROM with $m=1$ (i.e., low-order filtering) consistently produce
the most accurate results for large $N$ values, while a higher-order filter is
more effective for low $N$ values. In addition, L-ROM is sensitive to changes in
$N$ for all $m$ values. 
For $\rm Re_\tau =395$, L-ROM and TR-ROM with $m=1$ yield the most accurate
results for all $N$ values, while EFR-ROM yields the most accurate results for
low $N$ values for $m=1$, and high $N$ values for $m=3,4$.

To facilitate the comparison of the three Reg-ROMs, in
\cref{table:accuracy-ranking-pred}, we rank them based on the lowest error
achieved for the $N$ and $m$ values investigated.  Specifically, for both
Reynolds numbers, we list the Reg-ROMs' rank, the lowest $\euv$, the
corresponding $\euu$,  the ROM dimension $N$ and the filter order $m$ for which
the lowest $\euv$ is achieved.  The results in
\cref{table:accuracy-ranking-pred} yield the following conclusions:
\begin{table}[!ht]
\caption{
    Reg-ROM accuracy ranking in the predictive regime for $\rm Re_\tau = 180$
    (top rows) and $395$ (bottom rows).  The following parameters are listed: 
    Reg-ROMs' rank, the lowest $\euv$ with $\delta_\trecon$, the corresponding
    $\euu$, the ROM dimension $N$ and the filter order $m$ for which the lowest
    $\euv$ with $\delta_\trecon$ is achieved.
    \label{table:accuracy-ranking-pred} 
}
\centering
\begin{tabular}{|l|l|l|l|}
\hline
$\rm Re_\tau=180$  & L-ROM & EFR-ROM & TR-ROM \\ \hline
Rank & 2 & 3       & {1}      \\ \hline
$\euv$ & $\approx 22\%$    & $\approx 25\%$       & {$\approx 23\%$}      \\ \hline
$\euu$ & $\approx 23\%$    & $\approx 26\%$       & {$\approx 9\%$}      \\ \hline
N & $10$    & $80$ & $70 $      \\ \hline
Filter order $m$ &  1   & 1       & 1     \\ \hline
$\rm Re_\tau=395$ & L-ROM & EFR-ROM & TR-ROM \\ \hline
Rank & 2     & 3       & {1}      \\ \hline
$\euv$ & $\approx 39\%$    & $\approx 45\%$       & {$\approx 32\%$}      \\ \hline
$\euu$ & $\approx 39\%$    & $\approx 36\%$       & {$\approx 30\%$}      \\ \hline
N & $50$    & $40$ & $100$      \\ \hline
Filter order $m$ & $1$   & $1$       & $1$     \\ \hline
\end{tabular}
\end{table}

For $\rm Re_\tau=180$, TR-ROM is the most accurate model with $\euv \approx 23\%$ and $\euu \approx 9\%$, followed by L-ROM and EFR-ROM. 
For $\rm Re_\tau=395$, 
TR-ROM is still the most accurate model with $\euv \approx 32\%$ and $\euu \approx 30\%$, followed by EFR-ROM and L-ROM. 
In addition, the results in \cref{table:accuracy-ranking-pred} also show that $m=1$
(i.e., low-order filtering) yields the most accurate results.  Moreover, these
results show that TR-ROM requires large $N$ values to achieve its best accuracy,
whereas L-ROM yields best accuracy with small $N$. For EFR-ROM, large $N$ is
required for $\rm Re_\tau=180$, and small $N$ is required for $\rm Re_\tau =395$.

\begin{figure}[!ht]
    \centering
    \begin{subfigure}{0.49\columnwidth}
    \caption{L-ROM, $\rm Re_\tau=180$}
    \label{subfig:}
    \includegraphics[width=1\textwidth]{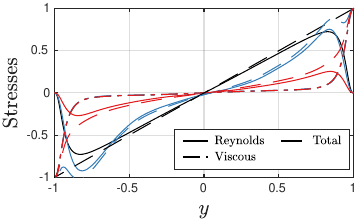}
    \end{subfigure}
    \begin{subfigure}{0.49\columnwidth}
    \caption{L-ROM, $\rm Re_\tau=395$}
    \label{subfig:}
    \includegraphics[width=1\textwidth]{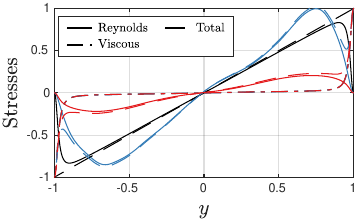}
    \end{subfigure}\\
    \begin{subfigure}{0.49\columnwidth}
    \caption{EFR-ROM, $\rm Re_\tau=180$}
    \label{subfig:}
    \includegraphics[width=1\textwidth]{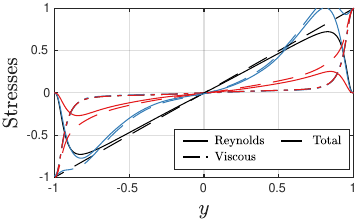}
    \end{subfigure}
    \begin{subfigure}{0.49\columnwidth}
    \caption{EFR-ROM, $\rm Re_\tau=395$}
    \label{subfig:}
    \includegraphics[width=1\textwidth]{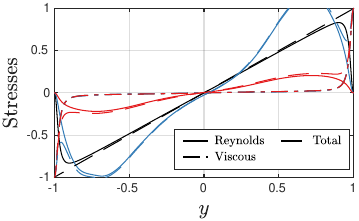}
    \end{subfigure}\\
    \begin{subfigure}{0.49\columnwidth}
    \caption{TR-ROM, $\rm Re_\tau=180$}
    \label{subfig:}
    \includegraphics[width=1\textwidth]{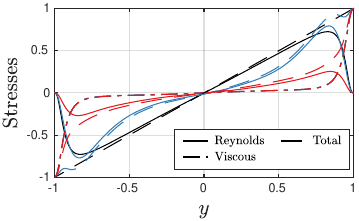}
    \end{subfigure}
    \begin{subfigure}{0.49\columnwidth}
    \caption{TR-ROM, $\rm Re_\tau=395$}
    \label{subfig:}
    \includegraphics[width=1\textwidth]{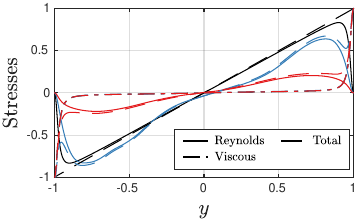}
    \end{subfigure}
    \caption{Comparison of the total, viscous, and Reynolds shear stresses of
    the FOM (black), the ROM projection (red), and the optimal Reg-ROMs, 
    listed in Table~\ref{table:accuracy-ranking-pred} (blue) in the \textit{predictive
    regime} for $\rm Re_\tau =180$ (left) and $\rm Re_\tau=395$ (right). The
    total shear stress is the sum of the viscous and the Reynolds shear stress.}
    \label{fig:pred_profiles}
\end{figure}

Similarly to the reproduction regime (Fig.~\ref{fig:recon_profiles}), in
Fig.~\ref{fig:pred_profiles}, we compare the total, viscous, and the Reynolds
shear stresses of the optimal Reg-ROMs (listed in
Table~\ref{table:accuracy-ranking-pred}) along with the results of the FOM and
the ROM projection in the predictive regime for $\rm Re_\tau=180$ and $\rm
Re_\tau=395$. 

In terms of the viscous shear stress, the results of the three Reg-ROMs are in
good agreement with those of the FOM and the ROM projection for both $\rm
Re_\tau$.  In terms of the Reynolds shear stress, we find that TR-ROM yields the
most accurate results for both $\rm Re_\tau$ and the L-ROM is better than the
EFR-ROM for $\rm Re_\tau=180$. For $\rm Re_\tau=395$, L-ROM and EFR-ROM perform
similarly.  In addition, from the results of the ROM projection, we find that $N=100$
POD basis functions are insufficient to reconstruct the Reynolds shear stress accurately.
Finally, as a results of the discrepancy in the Reynolds shear stress, the total
shear stress is not linear in all three Reg-ROMs for both $\rm Re_\tau$.

Next, we 
discuss the potential issue\sout{s} that preclude Reg-ROMs from getting more accurate 
Reynolds shear stress. 

From Figs.~\ref{fig:Re180_recon_conv}--\ref{fig:Re395_pred_conv}, we found that
$\euv$ of the ROM projection has at least $55\%$ errors for $\rm Re_\tau=180$,
and $70\%$ errors for $\rm Re_\tau=395$. 
    This indicates the poor approximation capability of the reduced basis functions for
    the Reynolds shear stress $\langle u'v' \rangle$.
With this approximation error, it is not surprising to find that $\euv$ of
Reg-ROMs is not below $10\%$. In order to improve the accuracy, one either needs
to increase the number of ROM basis functions or consider a ROM basis 
that is
designed for the Reynolds stress approximation.

\subsection{Sensitivity Study} 
    \label{section:sensitivity}

In this section, we perform sensitivity studies for the three Reg-ROMs: In
Section \ref{section:sensitivity_time_interval}, we present a sensitivity study
of the optimal parameter $(\delta,~\chi)_\trecon$ for each $N$ and $m$ 
in the predictive regime. 
In Section \ref{section:sensitivity_chi}, we present a sensitivity study of the
relative $\ell^2$ error $\euv$ of EFR-ROM and TR-ROM with respect to the
relaxation parameter, $\chi$.  In Section \ref{section:sensitivity_delta}, we
present a sensitivity study of the the relative $\ell^2$ error $\euv$ of the
three Reg-ROMs with respect to the filter radius, $\delta$.  
We refer
to the dissertation \cite{tsai2023pmordevelopment} for 
    the optimal filter radius $\delta_\trecon$'s sensitivity with respect to the filter order $m$ for the three Reg-ROMs.

\subsubsection{
Reg-ROM parameter sensitivity in the predictive regime}
    \label{section:sensitivity_time_interval}
    
In Section \ref{section:predictive-regime}, for given $N$ and $m$, we
investigated the accuracy of the three Reg-ROMs using
$(\delta,~\chi)_{\trecon}$ in the \textit{predictive regime}, where
$(\delta,~\chi)_\trecon$ are the $\delta$ and $\chi$ values that are optimal in
the reproduction regime. 
In this section, we extend that study and investigate
the robustness of the optimal parameter $(\delta,~\chi)_\trecon$ for each $N$
value with $m=1$ by comparing $(\delta,~\chi)_\trecon$ to
$(\delta,~\chi)_\tpred$, that is, the  $\delta$ and $\chi$ values that are
optimal in the predictive regime. 
We present $m=1$ results because this is the filter
order that yields the best models compared to other $m$ values. For the
sensitivity results of other $m$ values, see~\cite{tsai2023pmordevelopment}.
In order to
find 
$(\delta,~\chi)_\tpred$, we consider the same parameter sets for
$\delta$ and $\chi$ as those used in the reproduction regime (Section
\ref{section:reproduction-regime}).  An additional $500$ CTUs FOM simulations
are performed for both $\rm Re_\tau$ in order to compute $\euu$ and $\euv$. 

\begin{figure}[!ht]
    \centering
    \begin{subfigure}{0.49\columnwidth}
    \caption{L-ROM at $\rm Re_\tau = 180$}
    \label{subfig:re180_optdis_lrom}
    \includegraphics[width=1\textwidth]{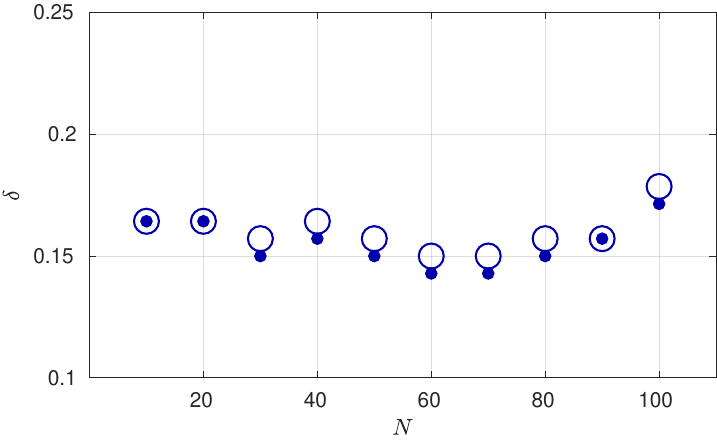}
    \end{subfigure}
    \begin{subfigure}{0.49\columnwidth}
    \caption{L-ROM at $\rm Re_\tau = 395$}
    \label{subfig:re395_optdis_lrom}
    \includegraphics[width=1\textwidth]{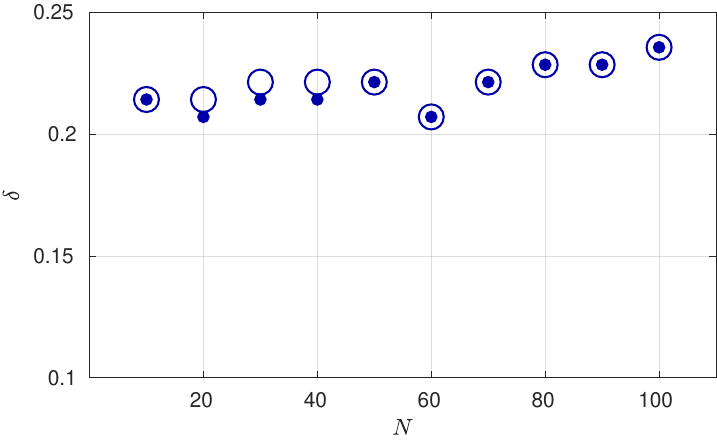}
    \end{subfigure}\\
    \begin{subfigure}{0.49\columnwidth}
    \caption{EFR-ROM at $\rm Re_\tau = 180$}
    \label{subfig:re180_optdis_efrrom}
    \includegraphics[width=1\textwidth]{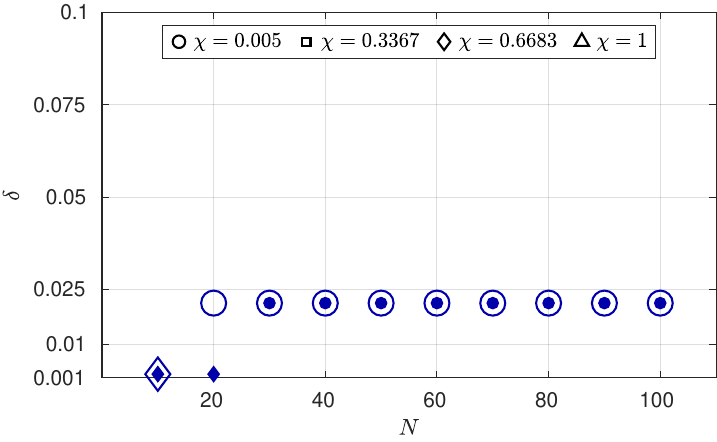}
    \end{subfigure}
    \begin{subfigure}{0.49\columnwidth}
    \caption{EFR-ROM at $\rm Re_\tau = 395$}
    \label{subfig:re395_optdis_efrrom}
    \includegraphics[width=1\textwidth]{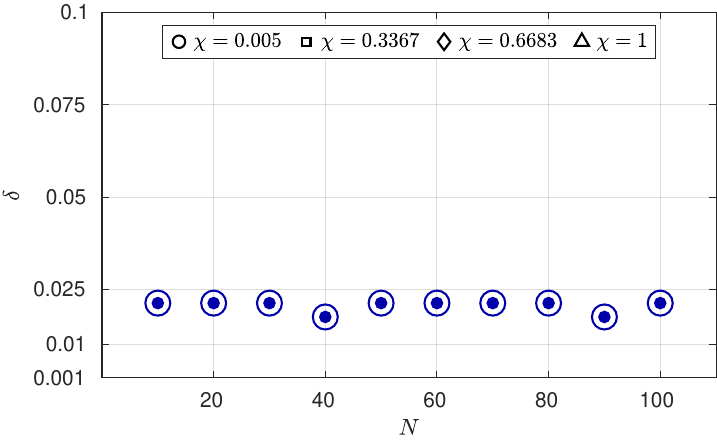}
    \end{subfigure}
    \begin{subfigure}{0.49\columnwidth}
    \caption{TR-ROM at $\rm Re_\tau = 180$}
    \label{subfig:re180_optdis_trrom}
    \includegraphics[width=1\textwidth]{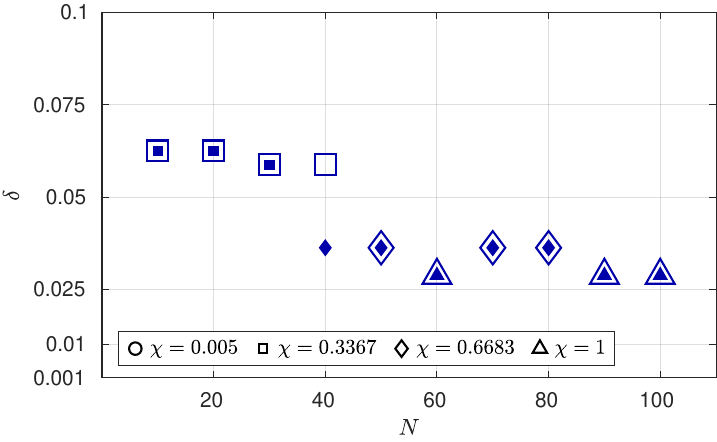}
    \end{subfigure}
    \begin{subfigure}{0.49\columnwidth}
    \caption{TR-ROM at $\rm Re_\tau = 395$}
    \label{subfig:re395_optdis_trrom}
    \includegraphics[width=1\textwidth]{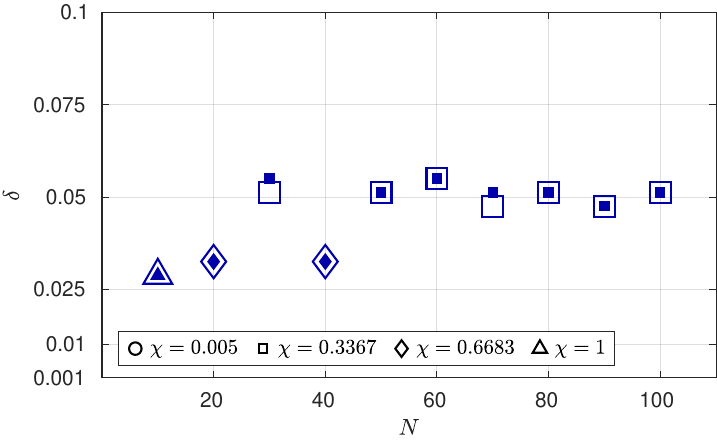}
    \end{subfigure}
    \caption{$(\delta,~\chi)_{\trecon}$ (smaller filled marker) and $(\delta,~\chi)_{\tpred}$ (larger empty marker)
    distributions of the three Reg-ROMs for $m=1$ and $N=10,\ldots,100$ at $\rm Re_\tau=180$ (left) and $\rm Re_\tau = 395$ (right).
    $(\delta,~\chi)_{\trecon}$ and $(\delta,~\chi)_{\tpred}$ are the optimal filter radius and relaxation parameter values found in the reproduction and predictive regimes, respectively.
    }
    \label{fig:optparam_distribution}
\end{figure}

In Fig.~\ref{fig:optparam_distribution}, we plot $\delta_\trecon$ and $\delta_\tpred$ for different $N$ values for the three Reg-ROMs at both $\rm Re_\tau$. 
$\delta_\trecon$ is displayed as a smaller filled marker, 
while $\delta_\tpred$ is displayed as a larger empty marker. 
We also use different markers to distinguish the relaxation parameter, $\chi$. 
We note that, for a given $N$, if $(\delta,~\chi)_\trecon$ is also optimal for the predictive regime, the two corresponding markers of different sizes will be on top of each other.

Figs.~\ref{subfig:re180_optdis_lrom} and \ref{subfig:re395_optdis_lrom} display L-ROM's $(\delta,~\chi)_\trecon$ and $(\delta,~\chi)_\tpred$ for each $N$ for $\rm Re_\tau=180$ and $\rm Re_\tau =395$, respectively. Note that the relaxation parameter, $\chi$, is 
not used in L-ROM (\ref{eqn:l-rom}) and, therefore, is not plotted. For $\rm Re_\tau =180$, we find that $\delta_\trecon$ is close to but not identical to $\delta_\tpred$ for most $N$ values. 
For $\rm Re_\tau =395$, we find that 
$\delta_\trecon$ is identical to $\delta_\tpred$ for most $N$ values, except for  $N=20,~30,~40$. 
In addition, for both $\rm Re_\tau$, we find that larger $N$ lead to  larger $\delta_\trecon$ and  $\delta_\tpred$.

Figs.~\ref{subfig:re180_optdis_efrrom} and \ref{subfig:re395_optdis_efrrom} display EFR-ROM's $(\delta,~\chi)_\trecon$ and $(\delta,~\chi)_\tpred$ for each $N$ for $\rm Re_\tau=180$ and $\rm Re_\tau =395$, respectively. 
For $\rm Re_\tau =180$, we find that  $(\delta,~\chi)_\trecon$ is optimal for almost all $N$ values except for $N=20$. 
For $\rm Re_\tau =395$, we find that $\delta_\trecon$ is identical to $\delta_\tpred$ for all $N$ values.
 
In addition, we find that the EFR-ROM's $\delta$ is less sensitive to $N$ compared to the L-ROM's $\delta$. 
Moreover, for the same $\chi_\trecon$, we find that $\delta_\trecon$ and $\delta_\tpred$ are not sensitive to $\rm Re_\tau$. 
Finally, we find that $\chi=0.005$ works best for most $N$ except for $N=10,~20$ at $\rm Re_\tau=180$.

Figs.~\ref{subfig:re180_optdis_trrom} and \ref{subfig:re395_optdis_trrom} display TR-ROM's $(\delta,~\chi)_\trecon$ and $(\delta,~\chi)_\tpred$ for each $N$ for $\rm Re_\tau=180$ and $\rm Re_\tau =395$, respectively. 
For $\rm Re_\tau =180$, we find that $\delta_\trecon$ is robust for all $N$ values except for $N=40$. 
For $\rm Re_\tau =395$, we again find that $\delta_\trecon$ is robust for all $N$ values except for $N=30$ and $N=70$.
Compared to the L-ROM and EFR-ROM results, we find that $(\delta,~\chi)_\trecon$ and  $(\delta,~\chi)_\tpred$ are in general more sensitive to $N$ and $\rm Re_\tau$ in TR-ROM.

In summary, across all three Reg-ROMs with $m=1$, we find that $\delta_\trecon$ is optimal for most $N$ values for both $\rm Re_\tau$, 
demonstrating the predictive capabilities of the three Reg-ROMs and their associated parameters.

\subsubsection{Reg-ROM sensitivity 
to the relaxation parameter $\chi$}
    \label{section:sensitivity_chi}

In this section, we study the EFR-ROM and TR-ROM sensitivity to the relaxation parameter, $\chi$. 
To this end, we consider $\euv$ as the metric. For each $N$ and $m$, we 
investigate how $\euv$ is affected by $\chi$, and what $\chi$ values yield the 
lowest $\euv$ values.

We consider four $\chi$ values, which are uniformly sampled in the interval $[\Delta t = 0.005,~1]$.
For each $N$ and $\chi$ values, we show $\euv$ with $\delta_\trecon$ and $m=1$.
We fix the filter order $m$ to be $1$ because this is the value that yields the best Reg-ROMs in the reproduction and predictive regimes (Sections \ref{section:reproduction-regime-summary} and \ref{section:predictive-regime-summary}). 
\begin{figure}[!ht]
    \begin{subfigure}{0.49\columnwidth}
    \caption{EFR-ROM, $\rm Re_\tau = 180$}
    \label{subfig:sen_chi_180_efr}
    \includegraphics[width=1\textwidth]{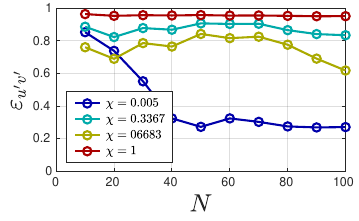}
    \end{subfigure}
    \begin{subfigure}{0.49\columnwidth}
    \caption{TR-ROM, $\rm Re_\tau = 180$}
    \label{subfig:sen_chi_180_tr}
    \includegraphics[width=1\textwidth]{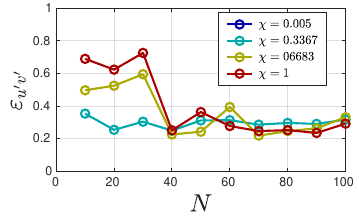}
    \end{subfigure}
    \begin{subfigure}{0.49\columnwidth}
    \caption{EFR-ROM, $\rm Re_\tau = 395$}
    \label{subfig:sen_chi_395_efr}
    \includegraphics[width=1\textwidth]{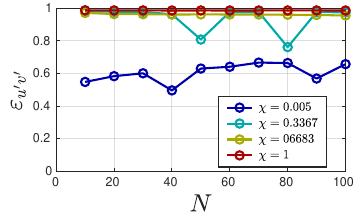}
    \end{subfigure}
    \begin{subfigure}{0.49\columnwidth}
    \caption{TR-ROM, $\rm Re_\tau = 395$}
    \label{subfig:sen_chi_395_tr}
    \includegraphics[width=1\textwidth]{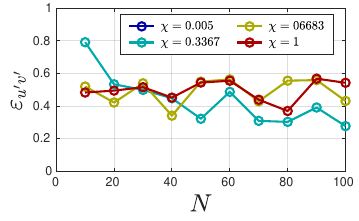}
    \end{subfigure}
    \caption{The relative error $\euv$ (\ref{equation:e-reystress}) of EFR-ROM
    and TR-ROM with respect to $N$ 
    and $\chi$ 
    for $\rm Re_\tau =180$ (\ref{subfig:sen_chi_180_efr} and
    \ref{subfig:sen_chi_180_tr}) and $\rm Re_\tau =395$
    (\ref{subfig:sen_chi_395_efr} and \ref{subfig:sen_chi_395_tr}) 
    with $m=1$ and $\delta_\trecon$,  
    where $\delta_\trecon$ is the optimal $\delta$ value found in the
    reproduction regime.}
    \label{fig:sen_wrt_chi}
\end{figure}

Fig.~\ref{subfig:sen_chi_180_efr} displays the EFR-ROM results at $\rm
Re_\tau=180$ for each $N$ and for four $\chi$ values. We recall that, in EFR-ROM
(\ref{eqn:efr-rom}), $\chi$ represents the contribution from the filtered
solution at each time step.  We find that EFR-ROM with $\chi = \Delta t = 0.005$
yields the best results. For $\chi=1$, EFR-ROM is too dissipative and leads to
an error of around $100\%$ in $\langle u'v' \rangle$ for all $N$ values. 
Although 
the number of samples we consider for $\chi$ is 
limited due to the training time, it is interesting to see that $\chi = \Delta t$ yields the best EFR-ROM results, just as in the FOM case 
\cite{ervin2012numerical,strazzullo2022consistency}.

Fig.~\ref{subfig:sen_chi_180_tr} displays the TR-ROM results for $\rm
Re_\tau=180$ for each $N$ and for four $\chi$ values. We recall that, in TR-ROM
(\ref{eqn:tr-rom}), $\chi$ represents the amount of additional diffusion added
to the G-ROM (\ref{eq:gromu}).  This time, we find that TR-ROM with $\chi =
\Delta t = 0.005$ yields the worst results. Because the amount of diffusion
added to G-ROM is not able to stabilize it, the error $\euv$ for each $N$ is
more than $100\%$.  For the other three $\chi$ values, for $N\le30$, we find
that smaller $\chi$ values lead to better accuracy, and for $N \ge 40$, we find
that the error $\euv$ is similar. 

Fig.~\ref{subfig:sen_chi_395_efr} displays the EFR-ROM results for $\rm
Re_\tau=395$ for each $N$ and for four $\chi$ values. We again find that EFR-ROM
with $\chi = \Delta t = 0.005$ yields the best results. In contrast with the
results for $\rm Re_\tau=180$, we find that EFR-ROM for the other three $\chi$
values leads to an error of around $100\%$ in $\langle u'v' \rangle$ for all $N$ values.

Fig.~\ref{subfig:sen_chi_395_tr} displays the TR-ROM results of $\rm
Re_\tau=395$ for each $N$ and for four $\chi$ values. Again, we find that TR-ROM
for $\chi = \Delta t = 0.005$ yields the worst results and the error $\euv$ is
more than $100\%$ for each $N$. For $N \le 40$, we find that $\chi = 0.6683$
yields the best results for $N=20$ and $40$, and $\chi=1$ yields the best
results for $N=10$ and $N=30$. For $N\ge50$, we find that $\chi=0.3367$ yields
the best results.

In summary, we found that both EFR-ROM and TR-ROM are sensitive to the
relaxation parameter $\chi$. Furthermore, for EFR-ROM, we found that
$\chi=\Delta t = 0.005$ outperforms the other three values for almost all $N$.
For TR-ROM, we found that $\chi=0.3367$ outperforms the other values for all $N$
for $\rm Re_\tau=180$, and for large $N$ for $\rm Re_\tau =395$.
\subsubsection{Reg-ROM sensitivity 
to the filter radius $\delta$}
    \label{section:sensitivity_delta}
    
In this section, for the optimal parameters listed in Table~\ref{table:accuracy-ranking}, we study the Reg-ROM sensitivity 
to the filter radius $\delta$ for $\rm Re_\tau =180$ and $\rm Re_\tau=395$. 
Our goal is to analyze the impact of $\delta$ on the Reg-ROM performance 
and identify the 
$\delta$ values that yield the best results.
\begin{figure}[!ht]
    \begin{subfigure}{0.49\columnwidth}
    \caption{L-ROM, $\rm Re_\tau = 180$}
    \label{subfig:re180_lrom_wrt_delta}
    \includegraphics[width=1\textwidth]{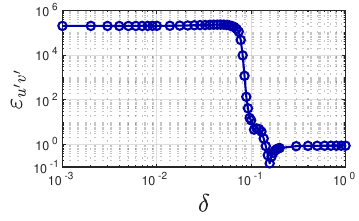}
    \end{subfigure}
    \begin{subfigure}{0.49\columnwidth}
    \caption{L-ROM, $\rm Re_\tau = 395$}
    \label{subfig:re395_lrom_wrt_delta}
    \includegraphics[width=1\textwidth]{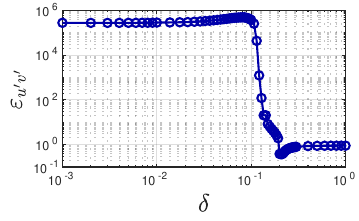}
    \end{subfigure}
    \caption{The relative error $\euv$ (\ref{equation:e-reystress}) of L-ROM 
    with respect to $\delta$ for $\rm Re_\tau=180$ (\ref{subfig:re180_lrom_wrt_delta}) and $\rm Re_\tau =395$ (\ref{subfig:re395_lrom_wrt_delta}).}
    \label{fig:sen_wrt_delta}
\end{figure}

Figures~\ref{subfig:re180_lrom_wrt_delta} and \ref{subfig:re395_lrom_wrt_delta}
display the L-ROM's $\euv$ behavior with respect to the filter radius, $\delta$,
for $\rm Re_\tau=180$ and $\rm Re_\tau =395$.  To discuss these results, we
divide the interval $[0.001,~1]$ into four subintervals: 
(i) For $\delta \in [0.001,~0.01]$, $\euv$ is high and does not change with respect to $\delta$. 
(ii) For $\delta \in [0.01,~0.1]$, $\euv$ increases as $\delta$ increases, and starts decreasing when $\delta$
approaches $0.01$. In addition, a much larger error drop is observed for $\rm
Re_\tau=180$ than for $\rm Re_\tau =395$. 
(iii) For $\delta \in [0.1,~0.2]$, $\euv$ decreases dramatically from $\mathcal{O}(10^5)$ down to $0.1$ and $0.37$ for $\rm Re_\tau = 180$ and $\rm Re_\tau=395$, respectively.
For both $\rm Re_\tau$, the optimal filter radius, $\delta_\trecon$, is obtained
in the interval $[0.1,~0.2]$. 
(iv) For $\delta \in [0.2,~1]$, $\euv$ increases as $\delta$ increases, and eventually 
plateaus at $\mathcal{O}(1)$ because L-ROM becomes too dissipative.

Figs.~\ref{subfig:re180_efrrom_wrt_delta} and
\ref{subfig:re395_efrrom_wrt_delta} display the EFR-ROM's $\euv$ behavior with
respect to the filter radius, $\delta$, for four $\chi$ values and for $\rm
Re_\tau=180$ and $\rm Re_\tau =395$. From the results, we can categorize
EFR-ROM's $\euv$ behavior into two types: 
\begin{figure}[!ht]
    \begin{subfigure}{0.49\columnwidth}
    \caption{EFR-ROM, $\rm Re_\tau = 180$}
    \includegraphics[width=1\textwidth]{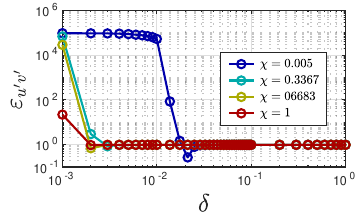}
    \label{subfig:re180_efrrom_wrt_delta}
    \end{subfigure}
    \begin{subfigure}{0.49\columnwidth}
    \caption{EFR-ROM, $\rm Re_\tau = 395$}
    \includegraphics[width=1\textwidth]{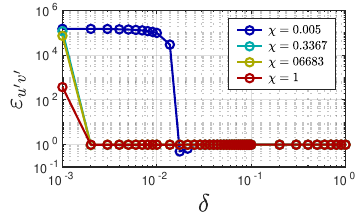}
    \label{subfig:re395_efrrom_wrt_delta}
    \end{subfigure}
    \caption{The relative error $\euv$ (\ref{equation:e-reystress}) of EFR-ROM 
    with respect to $\delta$ for $\rm Re_\tau=180$ (\ref{subfig:re180_efrrom_wrt_delta}) and $\rm Re_\tau =395$ (\ref{subfig:re395_efrrom_wrt_delta}).}
    \label{fig:sen_wrt_delta}
\end{figure}

For $\chi=\Delta t = 0.005$, the behavior of $\euv$ with respect to $\delta$ is
similar to that for L-ROM. We can divide the interval $[0.001,~1]$ into three
subintervals: (i) For $\delta \in [0.001,~0.01]$, $\euv$ decreases slightly as
$\delta$ increases.  (ii) For $\delta \in [0.01,~0.03]$, $\euv$ decreases
dramatically from $\mathcal{O}(10^5)$ down to $0.27$ and $0.5$ for $\rm Re_\tau
= 180$ and $\rm Re_\tau=395$, respectively.  For both $\rm Re_\tau$ values, the
optimal filter radius, $\delta_\trecon$, is obtained in the interval
$[0.01,~0.03]$.  (iii) For $\delta \in [0.03,~1]$, $\euv$ increases as $\delta$
increases, and eventually plateaus at $\mathcal{O}(1)$ because EFR-ROM becomes
too dissipative.

For the other three $\chi$ values, we can divide the interval $[0.001,~1]$ into
two subintervals: (i) For $\delta \in [0.001,~0.003]$, $\euv$ decreases
dramatically to $\mathcal{O}(1)$. However, the fact that there is no $\delta$
such that $\euv$ is below $1$ suggests that the optimal filter radius
$\delta_\trecon$ is either very sensitive, which requires more sampling points
in the interval $[0.001,~0.003]$, or it does not exist at all. (ii) For $\delta
\in [0.003,~1]$, $\euv$ is mostly $\mathcal{O}(1)$, suggesting that, for these
$\chi$ values, EFR-ROM is too dissipative regardless of the $\delta$ value. 

\begin{figure}[!ht]
    \begin{subfigure}{0.49\columnwidth}
    \caption{TR-ROM, $\rm Re_\tau = 180$}
    \label{subfig:re180_trrom_wrt_delta}
    \includegraphics[width=1\textwidth]{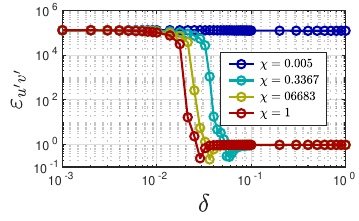}
    \end{subfigure}
    \begin{subfigure}{0.49\columnwidth}
    \caption{TR-ROM, $\rm Re_\tau = 395$}
    \label{subfig:re395_trrom_wrt_delta}
    \includegraphics[width=1\textwidth]{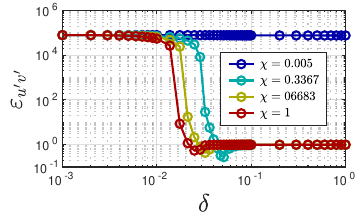}
    \end{subfigure}
     \caption{The relative error $\euv$ (\ref{equation:e-reystress}) of TR-ROM 
     with respect to $\delta$ for $\rm Re_\tau=180$ (\ref{subfig:re180_trrom_wrt_delta}) and $\rm Re_\tau =395$ (\ref{subfig:re395_trrom_wrt_delta}).}
    \label{fig:sen_wrt_delta}
\end{figure}

Figs.~\ref{subfig:re180_trrom_wrt_delta} and \ref{subfig:re395_trrom_wrt_delta}
display the TR-ROM's $\euv$ behavior with respect to the filter radius,
$\delta$, for four $\chi$ values and for $\rm Re_\tau=180$ and $\rm Re_\tau
=395$. 
From the results, we can categorize TR-ROM's $\euv$ behavior into two types:

For $\chi=\Delta t = 0.005$, there is no improvement in $\euv$ as $\delta$
increases. This suggests that $\chi$ is so small that, regardless of how large
the dissipative term $(\bu_r - \obu_r)$ (\ref{eqn:tr-rom-1}) (determined by
$\delta$) is, the total contribution $\chi (\bu_r - \obu_r)$ is too small to
stabilize G-ROM (\ref{eq:gromu}).

For the other three $\chi$ values, the TR-ROM's behavior of $\euv$ with respect
to $\delta$ is similar to that of L-ROM and EFR-ROM for $\chi=0.005$.  We can
divide the interval $[0.001,~1]$ into three subintervals:
(i) For $\delta \in [0.001,~0.01]$, $\euv$ decreases slightly as $\delta$ increases. 
(ii) For $\delta \in [0.01,~0.1]$, $\euv$ decreases dramatically from $\mathcal{O}(10^5)$ down to its optimal value. Furthermore, the smaller the $\chi$ value, the larger the $\delta$ value required to achieve its optimal $\euv$. This is expected because assuming that the total contribution $\chi (\bu_r - \obu_r)$ for optimal TR-ROM is fixed, larger $\chi$ will then require smaller $\delta$.
(iii) For $\delta \in [0.1,~1]$, $\euv$ increases as $\delta$ increases, and eventually plateaus at $\mathcal{O}(1)$ because TR-ROM becomes too dissipative.

In summary, we find that all three Reg-ROMs are sensitive to the filter radius,
$\delta$. For EFR-ROM and TR-ROM, $\delta_\trecon$ is affected by the relaxation
parameter $\chi$, and in the worst-case scenario, $\delta_\trecon$ might not
even exist.  In addition, we find that the optimal range for $\delta$ and the
effect of $\chi$ are similar for the two $\rm Re_\tau$ values.


\section{Conclusions and Future Work}
    \label{section:conclusions}

In this paper, we propose the time-relaxation ROM (TR-ROM), which is a novel
regularized ROM (Reg-ROM) for under-resolved turbulent flows. The TR-ROM employs
ROM spatial filtering to smooth out the flow velocity and eliminate the spurious
numerical oscillations displayed by the standard Galerkin ROM (G-ROM) (i.e., the
ROM that does not use any numerical stabilization).  
We emphasize that one novel feature of the TR-ROM, which distinguishes it from the
other Reg-ROM in current use (i.e., the Leray ROM (L-ROM) and the
evolve-filter-relax ROM (EFR-ROM)), is that it introduces different dissipation
for the large resolved scales and the small
resolved scales.  This is in stark contrast with the other two types of
Reg-ROMs, i.e, L-ROM and EFR-ROM, which use spatial filtering without
distinguishing between small and large resolved scales.

To assess the new TR-ROM, we compare it with the L-ROM and the EFR-ROM in the
numerical simulation of the turbulent channel flow at $Re_{\tau} = 180$ and
$Re_{\tau} = 395$ in both the reproduction and the predictive regimes.  The
spatial filtering in all three Reg-ROMs is performed using the first-order
ROM differential filter or the higher-order ROM algebraic filter.  We also
investigate the sensitivity of the Reg-ROMs with respect to the following
parameters: the time interval, the relaxation parameter, and the filter radius.
To our knowledge, this is the first numerical comparison of different Reg-ROMs
in the numerical simulation of turbulent flows.

Our numerical investigation yields the following conclusions: 
All three Reg-ROMs are dramatically more accurate than the classical G-ROM
without significantly increasing its computational cost.  
    In fact, with respect to several second-order turbulence statistics, the
    three Reg-ROMs' errors are much lower than 
    the projection error.
In addition, with the optimal parameters, the new TR-ROM yields more accurate results
than the L-ROM and the EFR-ROM in all tests. 
Our numerical investigation also shows that the HOAF with filter order $m=1$
yields the best results for most of the $N$ values. On the other hand, the HOAF
with $m>1$ works better for small $N$, i.e., $N\le 20$, at lower Reynolds number
$\rm Re_\tau=180$.

The sensitivity study shows 
that the optimal parameters trained in the
reproduction regime $(\delta,~\chi)_\trecon$ are also optimal in the predictive
regime for most of the $N$ values and for all three Reg-ROMs. 
Although 
all three Reg-ROMs are sensitive with respect to the
relaxation parameter, $\chi$, and filter radius, $\delta$, 
the optimal range for $\delta$ and the effect of $\chi$ are similar for the two
$\rm Re_\tau$ values. 

From the numerical investigation of the HOAF (Section
\ref{section:numerical-investigation-POD}), we found that the HOAF in the POD
setting is indeed a spatial filter, and has a similar behavior as in the SEM
setting, i.e., the higher-order filter (larger $m$) tends to damp the higher
index modes more, and has less impact on the lower index modes. 

The first steps in the numerical investigation of the new TR-ROM for the
turbulent channel flow have been encouraging.  There are, however, several other
research directions that should be pursued next.  For example, one can
investigate whether ROM approximate
deconvolution~\cite{sanfilippo2023approximate,xie2017approximate} can be
leveraged to further increase the localization of the dissipation mechanism in
the new TR-ROM.  One could also compare the approximate deconvolution approach
with the effect of increasing the order of the higher-order algebraic
filter~\eqref{eqn:hodf} in TR-ROM.  Finally, TR-ROM's numerical analysis, which
could yield new, robust parameter scalings~\cite{xie2018numerical}, should be
performed.

\section{Acknowledgments} 
This research used the Delta advanced computing and data resource which is
supported by the National Science Foundation (award OAC 2005572) and the State
of Illinois. The third author gratefully acknowledges NSF funding support
through grant DMS-2012253.

\appendix

\section{HOAF Investigation}
    \label{section:appendix}

In this section, we perform a theoretical and numerical investigation of the
HOAF~\eqref{eqn:hodf}, which, for clarity, we rewrite below:

{\em Given $\bu_\tr = \sum_{j=1}^N u_{\tr,j} \bphi_j$, 
find 
$\obu_\tr = 
\sum_{j=1}^N \ou_{\tr,j} \bphi_j$
such that  
\begin{eqnarray} 
    \left( \mathbbm{I} + \delta^{2m} A^{m} \right) \underline{\ou}_{\tr}
    = \uu_{\tr}, 	
    \label{equation:appendix-hodf}
\end{eqnarray}
where $\uu_{\tr}$ and $\underline{\ou}_{\tr}$ are the vectors of ROM basis
coefficients of $\bu_\tr$ and $\obu_\tr$, respectively, and $m \geq 2$ is an
integer.} 

We note that we consider only integers $m \geq 2$ because, as explained in
\cref{remark:hodf-m1}, the standard (low-order)
DF~\eqref{equation:df-linear-system} can be considered as a particular case of
HOAF with $m=1$.  Thus, to distinguish between the DF and the HOAF, in this
section we exclusively consider $m \geq 2$ in HOAF.

In addition, we note that the expansions for $\bu_\tr$ and $\obu_\tr$ do not
include the time-averaged velocity field, $\bphi_{0}$.  This is in contrast with
the expansion~\eqref{eq:romu}, which does include $\bphi_{0}$.  The reason for
not including $\bphi_{0}$ in our expansions is that, in our numerical
investigation, we decided not to filter the time-averaged velocity field,
$\bphi_{0}$, because this strategy was shown in~\cite{wells2017evolve} to yield
more accurate results. 

Moreover, we note that, in a spectral element setting,
in~\cite{mullen1999filtering} it was shown that
HOAF~\eqref{equation:appendix-hodf} is a spatial filter that attenuates the high
wavenumber components of the input field.  Furthermore, it was also shown that
the exponent $m$ in the HOAF~\eqref{equation:appendix-hodf} controls the percentage
of filtering at different wavenumbers: As $m$ increases, the amount of filtering
increases for the high wavenumber components of the input field, and decreases
for the low wavenumber components~\cite[Figure 1]{mullen1999filtering}.

In this section, we investigate whether the HOAF~\eqref{equation:appendix-hodf}
is a spatial filter.  We also study the role of the exponent $m$ in the
HOAF~\eqref{equation:appendix-hodf}.  We address this question first from a
theoretical point of view (Section~\ref{section:theoretical-investigation}), and
then from a computational point of view
(Section~\ref{section:numerical-investigation}).

\subsection{Theoretical Investigation}
    \label{section:theoretical-investigation}

In this section, we perform a theoretical investigation of HOAF~\eqref{equation:appendix-hodf}.
To this end, we first discuss the {\em weak form (variational formulation)} of HOAF, which is needed to construct the ROM discretization of HOAF.

We note that the weak form for DF is clear:
The DF linear systems~\eqref{equation:df-linear-system} results from the ROM discretization of the weak form of the Helmholtz equation~\eqref{equation:df-weak}.
Thus, because the DF can be interpreted as the HOAF~\eqref{equation:appendix-hodf} with $m=1$, we hope we can 
leverage that to find the HOAF weak form.

We emphasize that 
finding the HOAF weak form
is important.
Indeed, the connection between HOAF~\eqref{equation:appendix-hodf} and its conjectured 
weak form is used in the physical interpretation of the time relaxation term (Case 2 in Section~\ref{section:tr-rom}).
More importantly, the 
HOAF weak form 
can tell us whether HOAF is a spatial filter, just as in the DF case.

In this section, we 
discuss the HOAF weak form for
$m=2$. 
We believe this discussion can be naturally extended to higher $m$ values.

To discuss the HOAF weak form and the associated ROM formulation, 
we formally extend the mixed finite element discretization of the biharmonic equation proposed by Falk and Osborn~\cite{falk1980error}, which is based on earlier work by Ciarlet and Raviart~\cite{ciarlet1974mixed}, Glowinski~\cite{glowinski1973approximations}, and Mercier~\cite{mercier1974numerical}.

Because DF is HOAF with $m=1$, it is natural to associate the following PDE with HOAF:
{\em Given 
$\bu \in L^2(\Omega)$, find $\obu \in H_0^2(\Omega)$ such that}
\begin{eqnarray}
    && \delta^4 \Delta^2 \obu
    + \obu
    = \bu 
    \qquad \text{in } \Omega, 
    \label{equation:appendix-pde-strong} \\
    && \obu
    = \frac{\partial \obu}{\partial \bn}
    = 0
    \qquad \quad \text{on } \partial \Omega.
    \label{equation:appendix-pde-strong-bcs}
\end{eqnarray}

We note that, as pointed out in~\cite{ciarlet1974mixed}, the most straightforward finite element discretization of the biharmonic equation is to use a conforming method, in which the finite element space is a subspace of $H_0^2(\Omega) = \{\bv \in H^2(\Omega): \bv = \pp{\bv}{{\bm n}} = {\bm 0}~\text{on}~\partial \Omega\}$.
Constructing such subspaces, however, requires sophisticated finite elements.
Thus, to avoid these computational challenges, noncoforming finite elements have been proposed, in which the finite element space are subspaces of $H^1(\Omega)$.
One such strategy is the mixed finite element formulation of the biharmonic equation proposed in~\cite{falk1980error}, which utilizes standard finite element spaces in $H^1(\Omega)$.

We emphasize that, in our numerical investigation in Section~\ref{section:numerical-results}, the spectral element discretization used to construct the ROM basis function is posed in $H_0^1(\Omega)$.
Thus, we extend the mixed finite element strategy proposed in~\cite{falk1980error} to our ROM setting. 
To this end, we first define the following auxiliary variable:
\begin{eqnarray}
    \bw:= - \Delta \obu.
    \label{equation:appendix-1}
\end{eqnarray}
Next, we extend the weak form $(3.3)$ in~\cite{falk1980error} to the PDE we considered in~\eqref{equation:appendix-pde-strong}--\eqref{equation:appendix-pde-strong-bcs}:
{\em Given $\bu \in L^2(\Omega)$, find $(\bw,\obu) \in (H_0^1(\Omega), H_0^1(\Omega))$ such that}
\begin{eqnarray}
    && ( \bw , \bv_{\bw})
    - ( \nabla \obu , \nabla \bv_{\bw})
    = 0
    \qquad \qquad \quad \forall \, \bv_{\bw} \in H_{0}^{1}(\Omega) 
    \label{equation:appendix-weak-1} \\[0.3cm]
    && \delta^4 \, ( \nabla \bw , \nabla \bv_{\obu})
    + ( \obu , \bv_{\obu})
    = ( \bu , \bv_{\obu})
    \qquad \forall \, \bv_{\obu} \in H_{0}^{1}(\Omega) 
    \label{equation:appendix-weak-2}
\end{eqnarray}

\begin{remark}[ROM-FEM Weak Form]
    We note that, although the weak form~\eqref{equation:appendix-weak-1}--\eqref{equation:appendix-weak-2} is similar to the weak form $(3.3)$ in~\cite{falk1980error}, there are several significant differences.
    Probably the most important difference is that we use $H_0^1(\Omega)$ instead of $H^1(\Omega)$ to approximate the auxiliary variable $\bw$.
    Our choice is motivated by the spectral element space used in our numerical investigation in Section~\ref{section:numerical-results}, which is a subspace of $H_0^1(\Omega)$.
    \label{remark:falk}
\end{remark}

To find the ROM discretization of the weak formulations~\eqref{equation:appendix-weak-1} and \eqref{equation:appendix-weak-2}, we use the classical Galerkin method:
For~\eqref{equation:appendix-weak-1}, we consider the 
ROM expansions  
$\obu \approx \obu_\tr = \sum_{j=1}^N \ou_{\tr,j} \bphi_j$ and 
$\bw \approx \bw_\tr = \sum_{j=1}^N w_{\tr,j} \bphi_j$, 
and choose the test functions
$\bv_{\bw} := \bphi_i, i=1, \ldots, N$. 
This yields the following linear system: 
\begin{eqnarray}
    \mathbbm{I} \underline{w}_{\tr}
    = A \underline{\ou}_{\tr}, 
    \label{equation:appendix-7}
\end{eqnarray}
where $\underline{w}_{\tr}$ and $\underline{\ou}_{\tr}$ are the vector of ROM coefficients of $\bw_{\tr}$ and $\obu_{\tr}$, respectively, and $\mathbbm{I}$ and $A$ are the identity and ROM stiffness matrices, respectively.

To find the ROM discretization of the weak formulation~\eqref{equation:appendix-weak-2}, we use again the ROM expansions $\obu \approx \obu_\tr = \sum_{j=1}^N \ou_{\tr,j} \bphi_j$ and $\bw \approx \bw_\tr = \sum_{j=1}^N w_{\tr,j} \bphi_j$, and choose the test functions $\bv_{\obu} := \bphi_i, i=1, \ldots, N$.
This yields the following linear system: 
\begin{eqnarray}
    \delta^4 A \underline{w}_{\tr}
    + \mathbbm{I} \underline{\ou}_{\tr}
    = \mathbbm{I} \underline{u}_{\tr},
    \label{equation:appendix-8}
\end{eqnarray}
where $\underline{u}_{\tr}$ is the vector of ROM coefficients of $\bu_{\tr}$.
Plugging the formula for $\underline{w}_{\tr}$ in~\eqref{equation:appendix-7} in~\eqref{equation:appendix-8} yields 
\begin{eqnarray}
    \delta^4 A \left( A \underline{\ou}_{\tr} \right)
    + \mathbbm{I} \underline{\ou}_{\tr}
    = \mathbbm{I} \underline{u}_{\tr},
    \label{equation:appendix-9}
\end{eqnarray}
which is exactly the HOAF linear system~\eqref{equation:appendix-hodf} for $m=2$.

Thus, for $m=2$, we have shown that the ROM discretization of the weak forms~\eqref{equation:appendix-weak-1} and \eqref{equation:appendix-weak-2} is the HOAF linear system~\eqref{equation:appendix-hodf}.

We note that, in a spectral element setting, the PDE associated with HOAF~\eqref{equation:appendix-hodf} was discussed in~\cite{mullen1999filtering} (albeit with a missing minus sign).

\subsection{Numerical Investigation}
    \label{section:numerical-investigation}

In this section, we investigate whether the HOAF~\eqref{equation:appendix-hodf}
is a spatial filter by solving \eqref{equation:appendix-hodf} and examining if
the output $\obu$ is smoother compared to the input $\bu$. We also study the
role of the exponent $m$ in the HOAF.  In
sections~\ref{section:numerical-investigation-1DSEM}--\ref{section:numerical-investigation-2DSEM},
we investigate the HOAF in one- and two-dimensional spectral element (SEM)
settings. In section~\ref{section:numerical-investigation-POD}, we investigate
the HOAF in a POD setting.  We note that a similar numerical investigation was
carried out in the spectral element setting in~\cite{mullen1999filtering}, where
it was shown that HOAF~\eqref{equation:appendix-hodf} is a low-pass
filter~\cite[Figure 1]{mullen1999filtering}. 

\subsubsection{One-Dimensional SEM Setting}
    \label{section:numerical-investigation-1DSEM}

We consider the spatial domain $\Omega = [0,1]$ and construct the HOAF (\ref{equation:appendix-hodf}) 
using a 
SEM discretization that consists of a $64$-dimensional array of $7$th-order spectral elements.
Note that, in the SEM setting, instead of $A^m$ 
we use $(B^{-1}A)^m$ 
in (\ref{equation:appendix-hodf}), and 
$\uu$ and $\ubu$ are the SEM basis coefficients for the input (unfiltered) and output (filtered) functions, respectively.

\begin{figure}[!ht]
    \centering
        \includegraphics[width=0.9\linewidth]{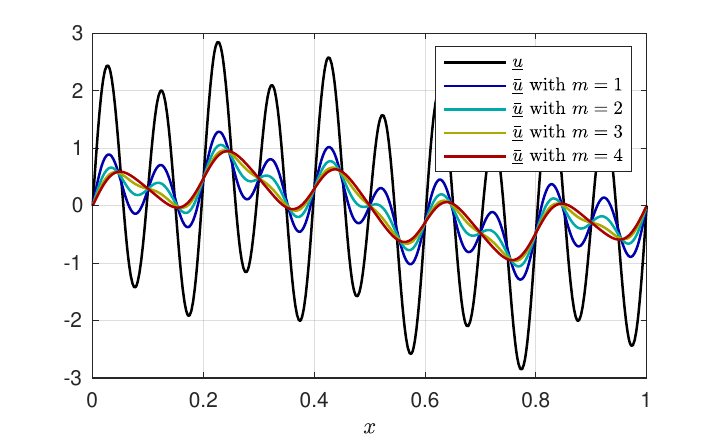}
    \caption{Comparison of the input function 
    $u(x)=0.5\sin(2\pi x)+ 0.5\sin(10 \pi x) $ $+ 2 \sin (20 \pi x)$ and four HOAF outputs $\ou$ that correspond to $m=1,~2,~3,~4$, and $\delta=0.025$.}
    \label{fig:hodf_1D_sem}
\end{figure}

We construct an input function $u$ 
that has both low and high wavenumber components:
\begin{equation}
    u(x)=0.5\sin(2\pi x) + 0.5\sin(10 \pi x) + 2 \sin (20 \pi x),
\end{equation}
and set the filter radius $\delta$ to be $0.025$.  We then consider four $m$
values and compare the corresponding HOAF output $\ou$ with the input function
$u$ in Fig.~\ref{fig:hodf_1D_sem}.

The $m=1$ results show that, although HOAF has attenuated the high wavenumber
($k=20$) component of the input function, that component is still visible.
However, as we increase the HOAF order (i.e., use $m=2, 3, 4$), the high
wavenumber component of the input function starts to significantly decrease.  In
fact, for $m=4$, the high wavenumber component is practically eliminated.  We
also note that the medium wavenumber ($k=10$) component of the input function is
affected differently by the different $m$ values: For $m=1$, the medium
wavenumber component is slightly attenuated.  
In contrast, for $m=4$, the medium wavenumber component remains unchanged.  
Finally, the low wavenumber ($k=2$) component of the input function is
practically unaffected by the HOAF order.  Thus, the results in
Fig.~\ref{fig:hodf_1D_sem} suggest that, as expected, using the classical DF
(i.e., HOAF with $m=1$) attenuates the medium and high wavenumber components of
the input function.  The higher-order HOAF (i.e., HOAF with $m=2,3, 4$)
attenuates more the high wavenumber component, and almost not at all the medium
wavenumber component. 

\subsubsection{Two-Dimensional SEM Setting}
\label{section:numerical-investigation-2DSEM}

\begin{figure}[!ht]
    \centering
    \begin{subfigure}{0.45\columnwidth}
    \caption{Input function $u$}
    \includegraphics[width=1\linewidth]{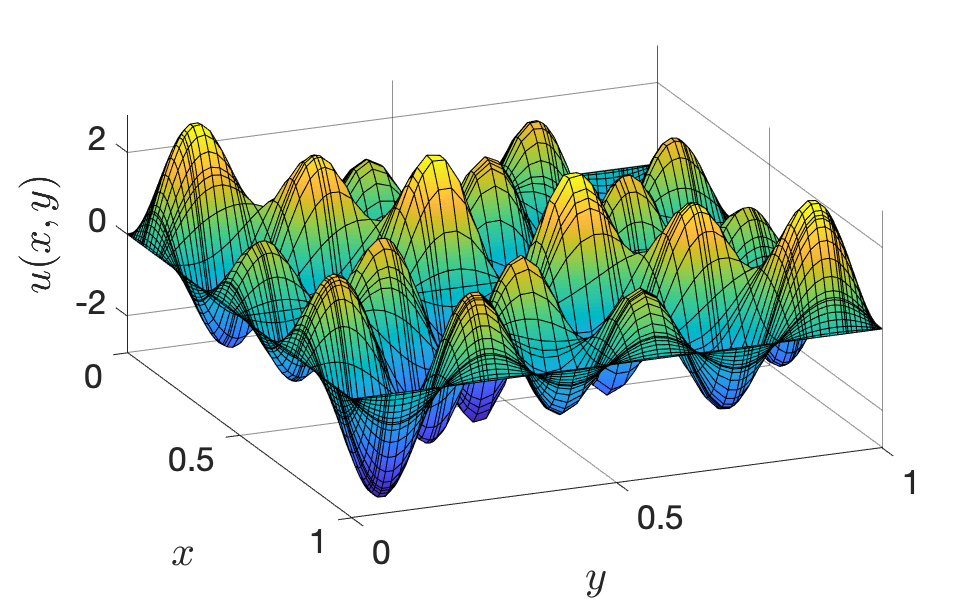}
    \end{subfigure}
    \begin{subfigure}{0.45\columnwidth}
    \caption{$\ou$ with $m=1$}
    \includegraphics[width=1\linewidth]{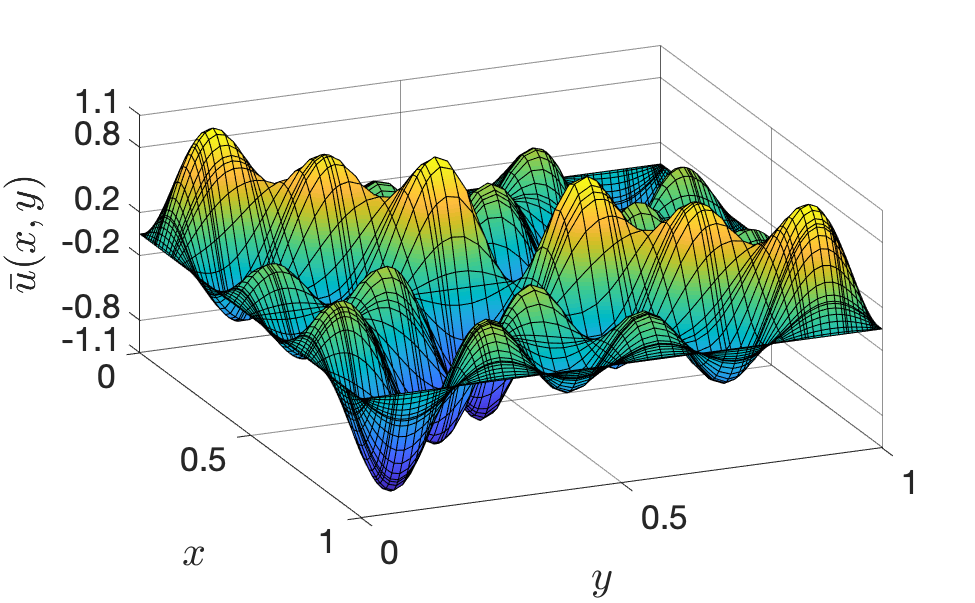}
    \end{subfigure}\\
    \begin{subfigure}{0.45\columnwidth}
    \caption{$\ou$ with $m=2$}
    \includegraphics[width=1\linewidth]{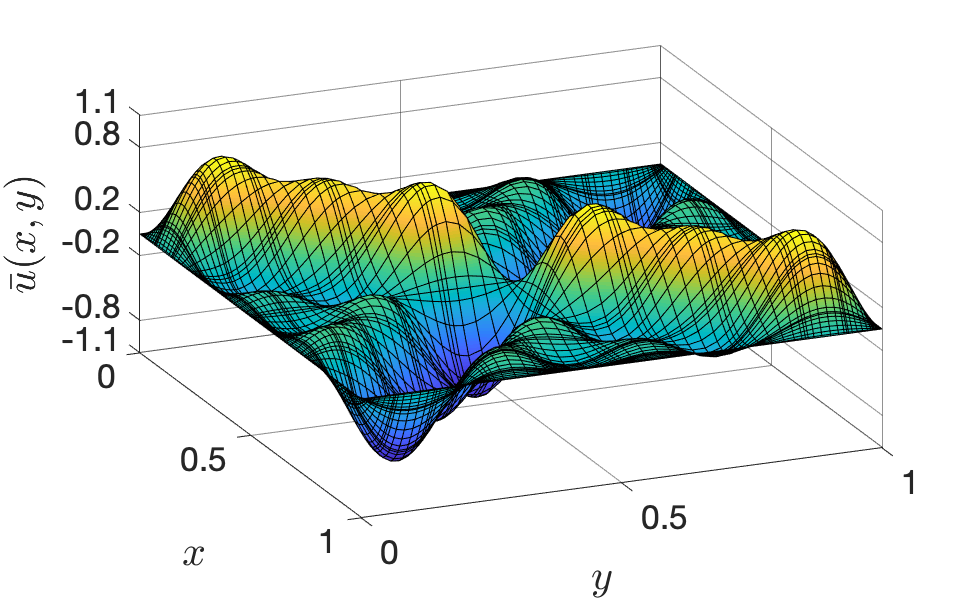}
    \end{subfigure}
    \begin{subfigure}{0.45\columnwidth}
    \caption{$\ou$ with $m=3$}
    \label{subfig:2Dsem-3}
    \includegraphics[width=1\linewidth]{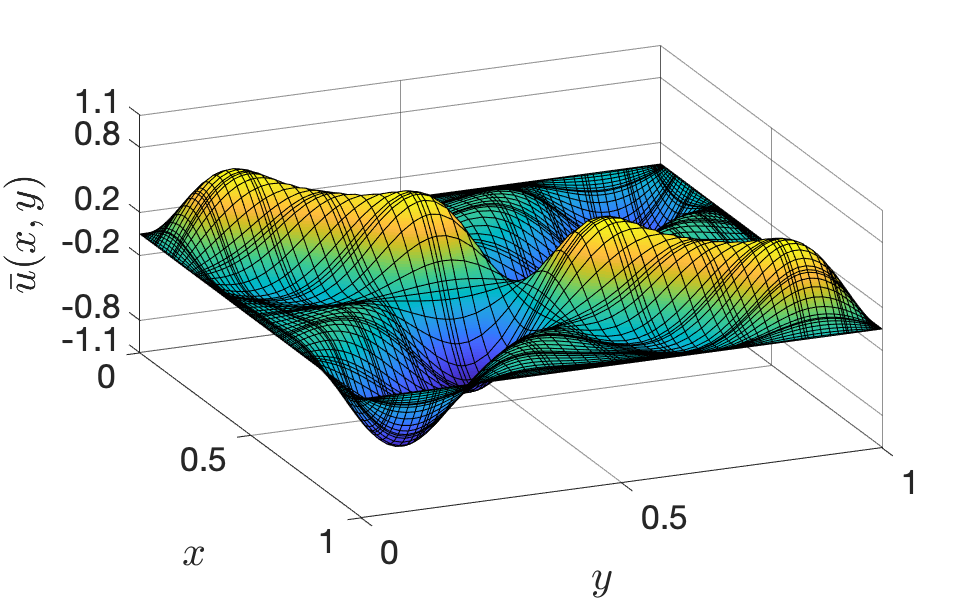}
    \end{subfigure}
    \caption{Comparison of the input function $u(x,y)=0.5\sin(2\pi x)\sin(2 \pi y) + \sin(4\pi x)\sin(4 \pi y)+ 2\sin(6\pi x)\sin(6 \pi y)$ and three HOAF outputs $\ou$ corresponding to $m=1,~2,~3$ with $\delta = 0.06$.}
    \label{fig:hodf_2D_sem}
\end{figure}

We consider the spatial domain $\Omega = [0,1]^2$ and construct the HOAF (\ref{equation:appendix-hodf}) 
using a SEM discretization that consists of a $12 \times 12$ array of
$7$th-order spectral elements.  We construct an input function $u$ that has both
low and high wavenumber components:
\begin{equation}
    u(x,y)=0.5\sin(2\pi x)\sin(2 \pi y) + \sin(4\pi x)\sin(4 \pi y)+ 2\sin(6\pi x)\sin(6 \pi y),
\end{equation}
and set the filter radius $\delta$ to be $0.06$. We then consider three $m$
values and compare the corresponding HOAF output $\ou$ with the input function
$u$ in Fig.~\ref{fig:hodf_2D_sem}.  
The $m=1$ results show that, although HOAF has attenuated the high wavenumber
($(k_1,k_2)=(6,6)$ component of the input function, that component is still
visible.  However, as we increase the HOAF order (i.e., use $m=2, 3$), the high
wavenumber component of the input function starts to significantly decrease.  We
also note that the medium wavenumber ($(k_1,k_2)=(4,4)$) component of the input
function is affected differently by the different $m$ values: For $m=1$, the
medium wavenumber component is slightly attenuated.  In contrast, for $m=3$, the
medium wavenumber component remains unchanged.  Finally, the low wavenumber
($(k_1,k_2)=(2,2)$) component of the input function is unaffected by the HOAF
order.  Thus, the results in Fig.~\ref{fig:hodf_2D_sem} again suggest that the
classical DF (i.e., HOAF with $m=1$) attenuates the medium and high wavenumber
components of the input function.  The higher-order HOAF (i.e., HOAF with
$m=2,3$) attenuates more the high wavenumber component, and almost not at all
the medium wavenumber component. 

\subsubsection{POD Setting}
\label{section:numerical-investigation-POD}

    We investigate the HOAF (\ref{equation:appendix-hodf}) in a POD setting.  In
    particular, we consider the POD basis functions 
    for the turbulent channel flow at
    $\rm Re_\tau=180$, as detailed in Section \ref{section:rom-computational
    setting}.  The HOAF is constructed 
    using $N=100$ 
    $L^2$ POD basis
    functions, and 
    filter radius $\delta = 0.08125$. The $\delta$ value is
    chosen such that the damping is neither too strong nor too weak.

    To understand how each POD basis function is affected by the HOAF with
    different $m$ values, we consider the following procedure:
    For a given $m$ value, and for each mode $i$, where $1\le i \le N$, we
    construct an input coefficient vector $\uu = [0,\ldots,0 ,1, 0,\ldots,0]^T$,
    where $1$ appears only in the $i$th 
    component (i.e., $u_i = 1$) and obtain the output
    vector $\ubu$ by solving the HOAF linear system (\ref{equation:appendix-hodf}.) The $i$th
    component of the output vector $\ubu$ (i.e., $\ou_i$) indicates the amount of
    damping caused by the HOAF in the $i$th POD basis function.
    %

\begin{figure}[!ht]
        \centering
        \includegraphics[width=0.8\linewidth]{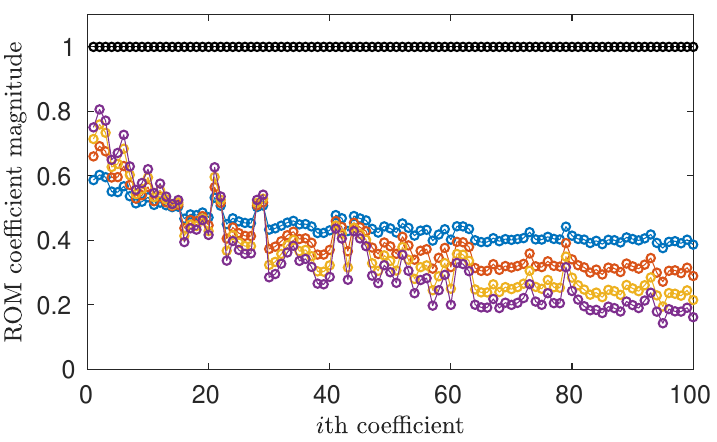}
        \caption{
            The unfiltered and filtered coefficient magnitude of the $i$th POD
            basis function, $u_i$ and $\ou_i$, with four filter order values
            $m=1,~2,~3,~4$ for $i=1,\ldots,N$.
            The HOAF (\ref{equation:appendix-hodf}) is formed using $N=100$
            $L^2$ POD basis functions and a filter radius of $\delta=0.08125$.
        }
        \label{fig:filter_coef_wrt_m}
    \end{figure}

    Fig.~\ref{fig:filter_coef_wrt_m} displays the unfiltered and the filtered
    coefficient magnitude of the $i$th POD basis function, $u_i$ (black) and
    $\ou_i$ (multi-colored), with four filter order values $m=1,~2,~3,~4$, for
    $i=1,\ldots,N$.
    Note that, for each mode $i$, 
    although the input vector $\uu$ 
    has only zero components 
    except 
    the $i$th 
    component (i.e., $u_i=1$), the output vector $\ubu$ could have
    a nonzero 
    $j$th 
    component (i.e., $\ou_j \neq 0$) for $j\neq i$, because
    the input vector $\uu$ is not necessary an eigenvector of the HOAF. 
    In Fig.~\ref{fig:filter_coef_wrt_m}, for each mode $i$, we  show only the
    behavior of $\ou_i$ because the purpose of this study is to investigate how
    each $i$th POD basis function is affected by the HOAF. In addition, we observe that, although
    $\ubu$ contains excitations of other modes, 
    these nonzero values are
    relatively small compared to the 
    $i$th 
    component. 
    %

    In Fig.~\ref{fig:filter_coef_wrt_m}, we see the expected behavior: 
    Increasing $m$ yields a sharper drop
    in the transfer function (${\bar u}_i/{u_i}$) for high mode numbers, $i$, with
    less suppression of the low mode numbers, which is consistent with an interpretation
    of the POD modes as ``Fourier'' modes.  Note that all the curves 
    intersect at 
    ${\bar u}_i/{u_i} \approx 1/2$ near $i=15$, which implies that modes with 
    $i < 15$ have characteristic length scales $ \lambda_i > \delta = 0.08125$,
    given that we expect 
    \begin{eqnarray*}
    \frac{{\bar u}_i}{u_i} & \sim & \frac{1}{1+(\delta/\lambda_i)^{2m}} \; .
    \end{eqnarray*}
    Thus, we see that, in this case, the HOAF provides a convenient means to 
    associate a length scale with each mode.

\begin{figure}[!ht]
        \centering
        \includegraphics[width=1\linewidth]{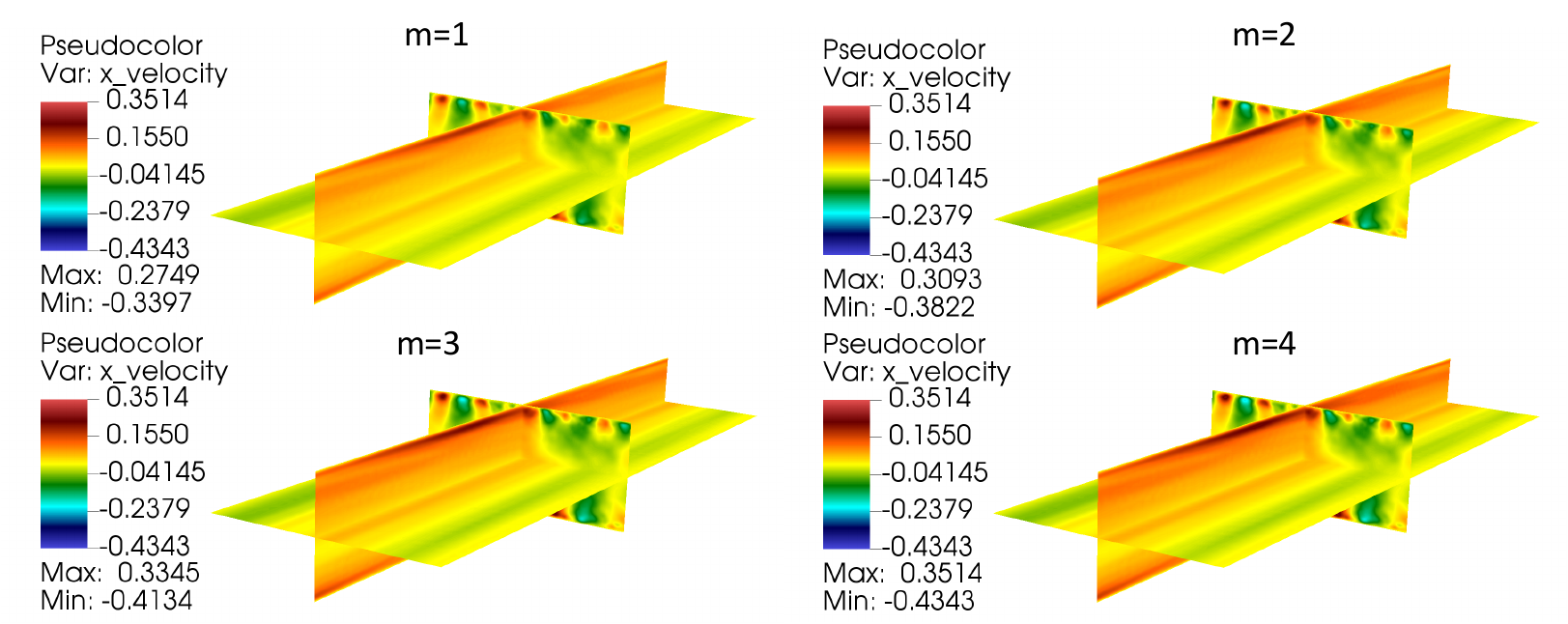}
        \caption{
            Comparison of streamwise component of the 1st filtered POD basis
            function for four filter order values $m=1,~2,~3,~4$. The HOAF is
            constructed with $N=100$ basis functions and filter radius
            $\delta=0.08125$.  The display range is fixed to be $[-0.4343,~0.3514]$
            and the maximum and minimum magnitudes are reported in the legend.
        }
        \label{fig:filter_pod_delta_0.08125_mode1}
    \end{figure}

\begin{figure}[!ht]
        \centering
        \includegraphics[width=1\linewidth]{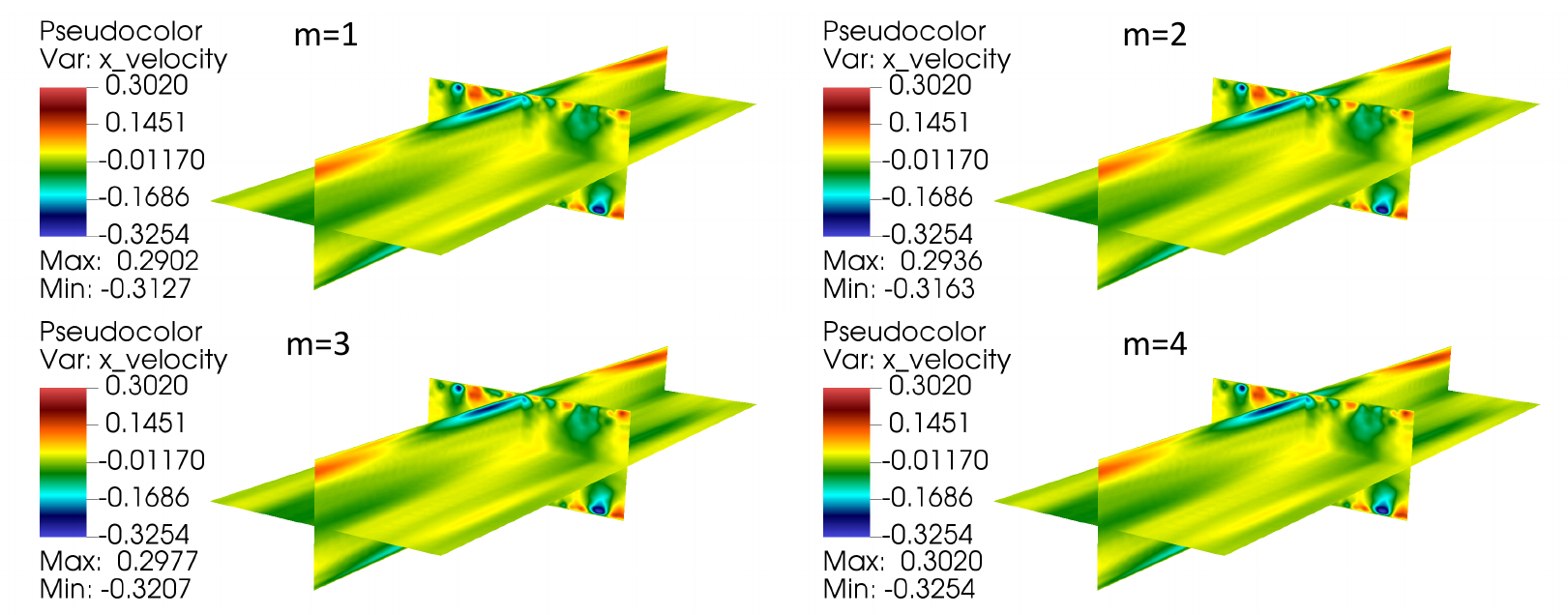}
        \caption{
            Comparison of streamwise component of the 15th filtered POD basis
            function for four filter order values $m=1,~2,~3,~4$. The HOAF is
            constructed with $N=100$ basis functions and filter radius
            $\delta=0.08125$.  The display range is fixed to be $[-0.3254,~0.3020]$
            and the maximum and minimum magnitudes are reported in the legend.
        }
        \label{fig:filter_pod_delta_0.08125_mode15}
    \end{figure}

\begin{figure}[!ht]
        \centering
        \includegraphics[width=1\linewidth]{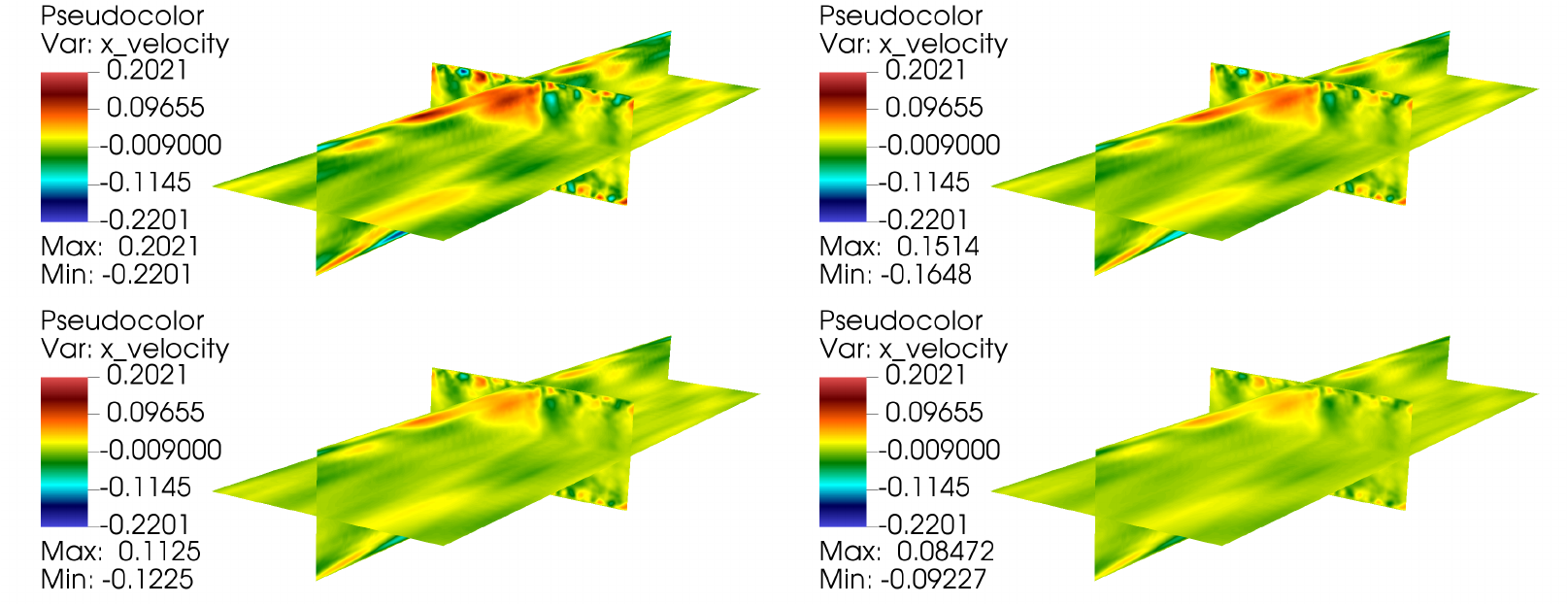}
        \caption{
            Comparison of streamwise component of the 100th filtered POD basis
            function for four filter order values $m=1,~2,~3,~4$. The HOAF is
            constructed with $N=100$ basis functions and filter radius
            $\delta=0.08125$.  The display range is fixed to be $[-0.2201,~0.2201]$
            and the maximum and minimum magnitudes are reported in the legend.
        }
        \label{fig:filter_pod_delta_0.08125_mode100}
    \end{figure}

    Figs.~\ref{fig:filter_pod_delta_0.08125_mode1}--\ref{fig:filter_pod_delta_0.08125_mode100}
    illustrate the physical-space effect of the HOAF with $m=1,~2,~3,~4$ and
    $\delta=0.08125$ for modes $i=1, 15$, and $100$.  It is clear that the
    higher-order filter (i.e., a larger $m$ value) has less damping in mode $i=1$
    and more damping in mode $i=100$. For mode $i=15$, all $m$ values yield
    similar damping.  Because the filtered basis function is simply the
    unfiltered basis function weighted by the coefficients given in
    Fig.~\ref{fig:filter_coef_wrt_m}, we expect the filtered basis function for
    different $m$ values 
    to differ only in the magnitude.
    Thus, in
    Figs.~\ref{fig:filter_pod_delta_0.08125_mode1}--\ref{fig:filter_pod_delta_0.08125_mode100},
    to compare the damping for different $m$ values, we compare the maximum and
    minimum magnitudes that are reported in the legends.

    Overall, these results show that the HOAF in the POD setting has a
    similar behavior as in the SEM setting, that is, larger $m$ values tend to damp the
    higher modes more and have less impact on the lower modes compared to smaller
    $m$ values. However, unlike in the SEM setting where the low-wave number modes
    are usually undamped regardless of the $m$ value, all POD modes are damped in
    the POD setting.

    In this numerical investigation we exclusively used $L^2$ POD
    basis functions, which is the most popular choice in reduced order modeling.  We note,
    however, that the $H^1_0$ POD basis functions could provide a clearer illustration of the
    HOAF spatial filtering capabilities.  Indeed, as highlighted in
    \cite{fick2018stabilized}, the $H^1_0$ POD basis functions are better in capturing the
    small-scale structures and distinguishing them from the large-scale structures
    in the solution compared to the $L^2$ POD basis functions.  Hence, we anticipate more
    distinct results in the filtered coefficients when the HOAF is constructed by
    using the $H^1_0$ POD basis functions.

\bibliographystyle{unsrt}
\bibliography{main}

\end{document}